\documentclass{mcrf}

\usepackage{amsmath,amsthm,amsfonts,amssymb,mathtools,mathrsfs,url}
  \usepackage{graphics} 
  \usepackage{epsfig}   
\usepackage{graphicx}   \usepackage{epstopdf}
\usepackage[colorlinks=true]{hyperref}
\hypersetup{urlcolor=blue, citecolor=red}

\def\currentvolume{X}
 \def\currentissue{X}
  \def\currentyear{200X}
   \def\currentmonth{XX}
    \def\ppages{X--XX}
     \def\DOI{10.3934/xx.xx.xx.xx}

\newtheorem{corollary}{Corollary}[section]
\newtheorem{definition}[corollary]{Definition}

\newtheorem{lemma}[corollary]{Lemma}
\newtheorem{proposition}[corollary]{Proposition}
\newtheorem{remark}[corollary]{Remark}
\newtheorem{theorem}[corollary]{Theorem}
\newcommand{\mylabel}[1]{\label{#1}
            \ifx\undefined\stillediting
            \else \fbox{$#1$}\fi }
\newcommand{\BE}{\begin{equation}}
\newcommand{\BEQ}[1]{\BE\mylabel{#1}}
\newcommand{\EEQ}{\end{equation}}
\newcommand{\rfb}[1]{\mbox{\rm
   (\ref{#1})}\ifx\undefined\stillediting\else:\fbox{$#1$}\fi}

\newfont{\Blackboard}{msbm10 scaled 1200}

\newfont{\sBlackboard}{msbm10 scaled 900}

\newfont{\roma}{cmr10 scaled 1200}
\renewcommand{\cline}{{\mathbb C}}
\newcommand{\fline}  {{\mathbb F}}
\newcommand{\nline}  {{\mathbb N}}
\newcommand{\rline}  {{\mathbb R}}
\newcommand{\tline}  {{\mathbb T}}

\renewcommand{\Re}{{\rm Re\,}}
\newcommand{\GGG}{{\bf G}}
\newcommand{\PPP}{{\bf P}}
\newcommand{\Dscr} {{\mathcal D}}
\newcommand{\Hscr} {{\mathcal H}}
\newcommand{\Lscr} {{\mathcal L}}
\newcommand{\Nscr} {{\mathcal N}}

\newcommand{\Zscr} {{\mathcal Z}}
\newcommand{\half} {{\frac{1}{2}}}
\newcommand{\mm}    {{\hbox{\hskip 0.5pt}}}
\newcommand{\m}     {{\hbox{\hskip 1pt}}}

\newcommand{\bluff} {{\hbox{\raise 15pt \hbox{\mm}}}}
\newcommand{\e}      {{\varepsilon}}
\renewcommand{\l}    {{\lambda}}
\renewcommand{\o}    {{\omega}}
\newcommand{\ka}     {{\kappa}} 
\newcommand{\FORALL} {{\hbox{$\hskip 6mm \forall \;$}}}
\newcommand{\kasten} {{\hfill\rule{2.5mm}{2.5mm}\hbox{\hskip 0.6mm}}}
\newcommand{\rarrow} {{\rightarrow}}

\newcommand{\jw}     {{i\omega}}
\newcommand{\dd}     {{\rm d\m}}

\newcommand{\LE}[1]    {{L^2([\mm 0,\infty),#1)}}
\newcommand{\LEloc}[1] {{L^2_{loc}([\mm 0,\infty),#1)}}
\newcommand{\Dom}[1]{{\mathcal D}(#1)}
\newcommand{\bigpar}[1]{\bigl( #1 \bigr)}
\newcommand{\CD}    {{C\& D}}
\newcommand{\bbm}[1]{\left[\begin{matrix} #1 \end{matrix}\right]}
\newcommand{\sbm}[1]{\left[\begin{smallmatrix} #1
            \end{smallmatrix}\right]}
\renewcommand{\theequation}{{\arabic{section}.\arabic{equation}}}

\font\fourteenrm = cmr10  scaled\magstep2
\font\fourteeni  = cmmi10 scaled\magstep2
\newcommand{\bsys} {{\hbox{\fourteenrm\char'006}}}
\newcommand{\bu}   {{\hbox{\fourteeni u}}}
\newcommand{\bv}   {{\hbox{\fourteeni v}}}
\newcommand{\by}   {{\hbox{\fourteeni y}}}
\newcommand{\bka}  {{\hbox{\fourteeni\char'024}}}
\title[Strong stabilization of passive systems]{Strong stabilization 
   of (almost) impedance passive systems by static output feedback}
\author[Ruth F. Curtain and George Weiss]{}
\email{gweiss@tauex.tau.ac.il}
\thanks{$^*$The second author is partially supported by the ETN
network ConFlex, funded by the European Union's Horizon 2020 research
and innovation programme under the Marie Sklodowska-Curie grant 
agreement No. 765579.}

\begin{document}
\maketitle

 \centerline{\scshape Ruth F. Curtain}
 \medskip {\footnotesize
 \centerline{Dept. of Mathematics}
 \centerline{University of Groningen}
 \centerline{9700 AV Groningen, The Netherlands}}

\medskip

 \centerline{\scshape George Weiss$^*$}
 \medskip {\footnotesize
 \centerline{School of Electrical Eng.}
 \centerline{Tel Aviv University}
 \centerline{Ramat Aviv 69978, Israel}}

\bigskip

\begin{abstract} The plant to be stabilized is a
system node $\Sigma$ with generating triple $(A,B,C)$ and transfer
function $\GGG$, where $A$ generates a contraction semigroup on the
Hilbert space $X$. The control and observation operators $B$ and $C$
may be unbounded and they are not assumed to be admissible. The
crucial assumption is that there exists a bounded operator $E$ such
that, if we replace $\GGG(s)$ by $\GGG(s)+E$, the new system
$\Sigma_E$ becomes impedance passive. An easier case is when $\GGG$ is
already impedance passive and a special case is when \mm $\Sigma$ has
colocated sensors and actuators. Such systems include many wave, beam
and heat equations with sensors and actuators on the boundary. It has
been shown for many particular cases that the feedback $u=-\ka y+v$,
where $u$ is the input of the plant and $\ka>0$, stabilizes $\Sigma$,
strongly or even exponentially. Here, $y$ is the output of \m $\Sigma$
and $v$ is the new input. Our main result is that if for some $E\in
\Lscr(U)$, $\Sigma_E$ is impedance passive, and \m $\Sigma$ is
approximately observable or approximately controllable in infinite
time, then for sufficiently small $\ka$ the closed-loop system is
weakly stable. If, moreover, $\sigma(A)\cap i\rline$ is countable,
then the closed-loop semigroup and its dual are both strongly
stable.
\end{abstract}

\section{The main results} 

\ \ \ This paper is a continuation and extension of our paper \cite
{CuWeABB*}, in the sense that we consider a more general class of
linear infinite-dimensional systems (system nodes) under more general
assumptions, as described below, and we are interested in stabilizing
these system nodes by static output feedback. The observability
assumption imposed here is weaker than in \cite{CuWeABB*}, and
correspondingly, the stability properties of the closed-loop system
will also be weaker. We shall refer often to background and results
from \cite{CuWeABB*}, but the concepts and techniques needed to derive
the results in this paper are quite different from those in
\cite{CuWeABB*}. 

To specify our terminology and notation, we recall that a {\it system
node} $\Sigma$ with input space $U$, state space $X$ and output space
$Y$ (all Hilbert spaces) is determined by its generating triple $(A,B,
C)$ and its transfer function $\GGG$, where the operator $A:\Dscr(A)
\rarrow\m X$ is the generator of a strongly continuous semigroup of
operators \m $\tline$ on $X$ and the possibly unbounded operators $B$
and $C$ are such that $C:\Dscr(A)\rarrow\m Y$ and $B^*:\Dscr(A^*)
\rarrow\m U$. There are no well-posedness assumptions for a system
node. Nevertheless, for input functions $u\in C^2([0,\infty),U)$ and
for initial states $z_0\in X$ that satisfy $Az_0+Bu(0)\in X$, there
exists a unique state trajectory $z\in C^1([0,\infty),X)$ satisfying
$\dot z =Az+Bu$. The corresponding $Y$-valued output function $y$ is
defined by $y=C[z-(\beta I-A)^{-1}u]+\GGG(\beta)u$ (where $\beta\in
\rho(A)$) and it is continuous. The space $\Dscr_0$ of all the pairs
$(z_0,u)$ which are as described above is dense in $X\times\LE U$. For
$(z_0,u)\in\Dscr_0$, the functions $z,y$ can also be expressed in
terms of their Laplace transforms:
$$ \begin{aligned} \hat z(s) &= (sI-A)^{-1}[z_0 + B\hat u(s)],\\
   \hat y(s) &= C(sI-A)^{-1} z_0 + \GGG(s)\hat u(s), \end{aligned}$$
for all $s\in\cline$ with sufficiently large real part. For details on
system nodes we refer to Malinen {\it et al} \cite{MalStafWei}, Opmeer
\cite{Opm}, Staffans \cite{stafbook}, and a short introduction to
system nodes will be given in Section 2.

The system node $\Sigma$ is called {\it impedance passive} if \m $Y=U$
and for all $(z_0,u)$ in the space $\Dscr_0$ (defined a little
earlier) and for all $\tau>0$,
\BEQ{statepassive}
   \|z(\tau)\|^2-\|z_0\|^2\leq 2\int_0^\tau {\rm Re}\langle u(t),
   y(t) \rangle \dd t \m.
\EEQ
Necessary and sufficient conditions for passivity in terms of $A,B,C,
\GGG$ will be recalled in Section 3. The concept of impedance
passivity has been introduced in Willems \cite{Will} in the
finite-dimensional context and has been generalized to the
infinite-dimensional context (with modified terminology) by several
researchers, see Staffans \cite{StafPas,StafCol} and the references
therein.

In this paper we deal with systems (plants to be stabilized) that
satisfy a relaxed version of \rfb{statepassive}: there exists an
$E\in\Lscr(U)$ such that
\BEQ{relaxed}
   \|z(\tau)\|^2-\|z_0\|^2\leq 2\int_0^\tau {\rm Re}\langle u(t),y(t)
   \rangle\m\dd t + 2\int_0^\tau \langle Eu(t),u(t) \rangle\m\dd t \m,
\EEQ
for all $(z_0,u)\in\Dscr_0$ and all $\tau\geq 0$. Equivalently, if we
replace $\GGG$ by $\GGG+E$ (and keep $A,B,C$ unchanged), then we
obtain a modified system node $\Sigma_E$ which is impedance passive.
We call such system nodes {\it almost impedance passive}.

We remark that the existence of $E$ such that \rfb{relaxed} holds is
obviously equivalent to the existence of $c\geq 0$ such that
$$ \|z(\tau)\|^2-\|z_0\|^2\leq 2\int_0^\tau {\rm Re}\langle u(t),y(t)
   \rangle\m\dd t + 2c\int_0^\tau \|u(t)\|^2 \dd t \m.$$
This condition is simpler, but for certain arguments it is better to
refer to \rfb{relaxed}.

In \cite{CuWeABB*} we considered the class of well-posed linear
systems with colocated actuators and sensors (meaning that $C=B^*$)
and $A$ was essentially skew-adjoint. In Section 3 we revisit this
class but without the well-posedness assumption. We show that these
systems satisfy \rfb{relaxed} and moreover, for these systems we
determine the minimal $E$ so that \rfb{relaxed} holds. 

Although the $A$ operator of an almost impedance passive system node
always generates a contraction semigroup, it is not necessary that the
actuators and sensors are colocated. Moreover, a system node with
colocated actuators and sensors and a contraction semigroup need not
be impedance passive. We give concrete examples that illustrate these
facts in Sections 6, 7 and 8, where we investigate in detail the
conditions for almost impedance passivity for some classes of second
order systems.

In Section 4 we consider an impedance passive system node and we
examine the effect of the static output feedback $u=-\ka y+v$, where
$\ka>0$ and $v$ is the new input function, as shown in a block diagram
in Figure 1. We give conditions under which this feedback results in a
closed-loop system $\Sigma^\ka$ that is well-posed and has nice
stability properties. To state our main results, we recall some
stability concepts.

\begin{definition} {\rm
Let $\Sigma$ be a well-posed linear system with input space $U$, state
space $X$, output space $Y$, generating triple $(A,B,C)$ and transfer
function $\GGG$.
\begin{itemize}
\item[$\bullet$] $\Sigma$ is {\it input stable} if for any $v\in\LE U
   $, the state trajectory of \m $\Sigma_0$ corresponding to the
   initial state zero and the input function $v$ is bounded. This
   property is also known as {\it infinite-time admissibility} of $B$.
\item[$\bullet$] $\Sigma$ is {\it output stable} if for any $z_0\in
   X$, the output function of \m $\Sigma$ corresponding to the initial
   state $z_0$ and the input function zero is in $\LE Y$. This
   property is also known as {\it infinite-time admissibility} of $C$.
\item[$\bullet$] $\Sigma$ is {\it input-output stable} if for any $v
   \in\LE U$, the output function of \m $\Sigma$ corresponding to the
   initial state zero and the input function $v$ is in $\LE U$.
   Equivalently, $\GGG\in H^\infty(\Lscr(U))$, the space of bounded
   analytic $\Lscr(U)$-valued functions on the open right half-plane.
\item[$\bullet$] $\Sigma$ is {\it system stable} if it is input
   stable, output stable and input-output stable.
\end{itemize}}
\end{definition}

We shall also use the following standard stability concepts for a
strongly continuous semigroup $\tline$ on a Hilbert space $X$:
\begin{itemize}
\item[$\bullet$] {\it weak stability} means that $\langle\tline_t
   z_0,z_1\rangle\m\rarrow\m 0$ as $t\m\rarrow\m\infty$, for all
   $z_0,z_1\in X$,
\item[$\bullet$] {\it strong stability} means that $\tline_t z_0\m
   \rarrow\m 0$ as $t\m\rarrow\m\infty$, for all $z_0\in X$,
\item[$\bullet$] {\it exponential stability} means that there exist
   $M\geq 1$ and $\alpha>0$ such that \\ $\|\tline_t\|\leq M
   e^{-\alpha t}$ for all $t\geq 0$.
\end{itemize}

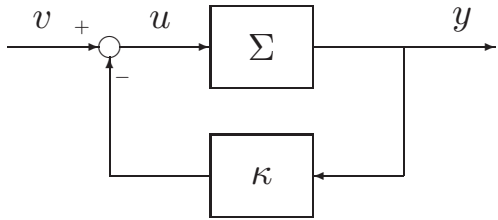
\begin{figure}[bth]
\unitlength 9mm
\begin{center}
\begin{picture}(6.5,3)
\put(1.0,0.4)  {\vector(0,1){1.75}}            
\put(1.2,1.95) {\makebox(0,0)[t]{$_-$}}        
\put(1.0,2.3)  {\circle{0.3}}                  
\put(-.5,2.3)  {\vector(1,0){1.35}}            
\put(0.6,2.5)  {\makebox(0,0)[b]{$_+$}}        
\put(1.15,2.3) {\vector(1,0){1.33}}            
\put(2.5,1.7)  {\framebox(1.5,1.2){$\bsys$}}   
\put(4.0,2.3)  {\vector(1,0){2.7}}             
\put(5.35,2.3) {\line(0,-1){1.9}}              
\put(5.35,0.4) {\vector(-1,0){1.34}}           
\put(2.5,-0.2) {\framebox(1.5,1.2){$\bka$}}    
\put(2.5,0.4)  {\line(-1,0){1.5}}              
\put(0.0,2.7)  {\makebox(0,0){$\bv$}}          
\put(1.75,2.7) {\makebox(0,0){$\bu$}}          
\put(6.2,2.7)  {\makebox(0,0){$\by$}}          
\end{picture}
\end{center}
\vskip -1mm
\caption{The open-loop system node $\Sigma$ with static output
feedback. If $\ka>0$ is sufficiently small, then this feedback results
in a closed-loop system $\Sigma^\ka$ that is well-posed and system
stable. Under suitable additional assumptions, the operator semigroup
of \m $\Sigma^\ka$ and its dual are strongly stable.
\hfill \m \hskip 90mm \m\hskip 37mm --------------------------}
\end{figure}

A consequence of our results in Section 4 is the following. 

\begin{theorem} \label{inpoutstab}
Let \m $\Sigma$ be a system node for which there exists $E\in\Lscr(U)$
satisfying \rfb{relaxed}. Then there exists $\ka_0>0$ (possibly $\ka_0
\mm=\infty$) such that for all $\ka\in(0,\ka_0)$, the feedback law
$u=-\ka y+v$ (where $u$ and $y$ are the input and the output of \m
$\Sigma$) leads to a well-posed closed-loop system \m $\Sigma^\ka$
that is system stable.

Moreover, the semigroup of \ $\Sigma^\ka$ is a contraction semigroup.
\end{theorem}

With the assumptions of the above theorem, denote $\Re E=\half(E+E^*)
$. If $c$ is the smallest number in $[0,\infty)$ such that $\Re E\leq
cI$, then we can take $\kappa_0=1/c$ (if $c=0$ then $\kappa_0=\infty
$). In particular, if \m $\Sigma$ is impedance passive then
$\kappa_0=\infty$.

In Section 5 we prove some results concerning weak and strong
stability under mild observability and controllability assumptions. We
say that $(A,C)$ (or the system node $\Sigma$) is {\it approximately
observable in infinite time} if for every $z_0\in\Dscr(A),$ $C\tline_t
z_0=0\ \forall t\geq 0$ implies $z_0=0$. We say that $(A,B)$ (or the
system node $\Sigma$) is {\it approximately controllable in infinite
time} if $(A^*,B^*)$ is approximately observable in infinite time. The
following theorem is a consequence of our results in Section 5.

\begin{theorem} \label{approx_obs}
Under the notation and assumptions of Theorem {\rm\ref{inpoutstab}},
if $(A,B,C)$ is either approximately controllable in infinite time or
approximately observable in infinite time, then for every $\ka\in(0,
\ka_0)$, the closed-loop semigroup $\tline^\ka$ is weakly stable. If,
in addition, the intersection of the spectrum $\sigma(A)$ with the
imaginary axis $i\rline$ is at most countable, then $\tline^\ka$ and
its dual $\tline^{\ka*}$ are strongly stable.
\end{theorem}

The proof of the strong stability part of this theorem uses a famous
result from Arendt and Batty \cite{AreBat} and Lyubich and Phong
\cite{LyPh} (see also Arendt {\it et al} \cite{ABHN}).

In the literature we can find various particular examples or classes
of passive systems stabilized by static output feedback. The result
most closely related to the last theorem is the main result of Batty
and Phong \cite{BatPho}. In that paper, it is assumed that $A$
generates a contraction semigroup, $B$ is bounded (i.e., it is in
$\Lscr(U,X)$), $C=B^*$ and the open-loop transfer function is $\GGG
(s)=B^*(sI-A)^{-1}B$. Similar assumptions have been made in a 
series of earlier papers: The first PDE examples fitting into this
framework were Bailey and Hubbard \cite{BaiHub}, Balakrishnan
\cite{Bal1,Bal2}, Russell \cite{Rus1}, Slemrod \cite{Slem2,Slem3}. In
this case, $\GGG$ is positive and (by a simple argument) this implies
that the feedback $u=-\ka y+v$ stabilizes in an input-output sense. Of
course, the most desirable type of stability is exponential stability
and for this, in the special case $A^*+A=0$, we need the system to be
exactly controllable (or equivalently, exactly observable). This is
the setup studied in Haraux \cite{Haraux}, Liu \cite{Liu97}, Lasiecka
and Triggiani \cite{LasTrig}, \cite{LasTrig3} and others.

However, early on it was realized that exponential stability is not
achievable with a bounded $B$ in the case that $U$ is
finite-dimensional and $A$ has infinitely many eigenvalues on the
imaginary axis. This is illustrated in Russell \cite{Rus1} with a PDE
model of an undamped string. Similar results for a beam can be found
in Slemrod \cite{Slem3}. More generally, it is known that if $U$ is
finite-dimensional, $B$ is bounded and $A$ has infinitely many
unstable eigenvalues, we can never achieve exponential stability, see
Gibson \cite{Gib}, Triggiani \cite{trigg2} or Curtain and Zwart
\cite[Theorem 5.2.6]{CurZwa95}. So for this class the best we can hope
for is strong stability. Early results giving sufficient conditions
under which $A-BB^*$ generates a weakly or strongly stable semigroup
using LaSalle's principle can be found in Slemrod \cite{Slem1}. These
were sharpened by Benchimol in \cite{Ben} who used the canonical
decomposition of contraction semigroups due to Sz\"okefalvi-Nagy and
Foias \cite{NagFoi}. He showed that if $A$ generates a contraction
semigroup and $B\in\Lscr(U,X)$, a sufficient condition for $A-BB^*$ to
generate a weakly stable semigroup is that
\BEQ{bench}
   \left\{ x\in X \ |\ B^*\tline^*_t x = 0 \m,\ \
   \|\tline_t x\|=\|x\|=\|\tline^*_tx\| 
   \FORALL t > 0 \right\} \m=\m \{ 0 \} \m.
\EEQ
In \rfb{bench}, $\tline^*_t$ may be interchanged with $\tline_t$. If,
in addition, $A$ has compact resolvents, then \rfb{bench} implies
strong stability. The above result for weak stability was also
obtained by Batty and Phong \cite{BatPho}. They improved the above
sufficient condition for strong stability, obtaining: if the spectrum
of $A$ has at most countably many points of intersection with the
imaginary axis, then $A-BB^*$ generates a strongly stable semigroup
if and only if \rfb{bench} holds. It is worthwhile noting that, while
the assumption that $B$ be bounded is restrictive, it does not exclude
PDEs with boundary control, see Slemrod \cite{Slem3}, You \cite{You}
and Chapter 9 of Oostveen \cite{Oos}.

As already mentioned, there are systems that are impedance passive,
but $C^*\neq B$. In \cite{LuoGuoMor98} Luo, Guo and Morgul gave
conditions for a class of PDEs to be {\em impedance energy preserving}
(i.e., \rfb{statepassive} holds with equality). The PDEs were in one
spatial dimension with control and observation on the boundary and $A$
was assumed to be skew-adjoint and to have compact resolvents. Using
LaSalle's principle they showed that under an observability
assumption, the feedback $u=-\ka y+v$ with $\ka>0$ produces a strongly
stable closed-loop system. Using a different approach Le Gorrec {\it
et al} \cite{GorZwaMas05} defined a large class of hyperbolic-like
systems in one spatial variable with boundary control and observation
that were impedance energy preserving. In Zwart {\it et al} \cite
{gang} it was shown, using a Lyapunov approach, that in the constant
coefficients case these systems could be stabilized by static output
feedback.

Most of the existing results on static output stabilization of PDE
systems assume a skew-adjoint semigroup generator and colocated
actuators and sensors. Theorems \ref{inpoutstab} and \ref{approx_obs}
include and generalize these results; neither a skew-adjoint generator
nor colocated actuators and sensors are needed.

{\bf Historical note.} This paper has been written over many years, up
to August 2007. The authors wanted to improve certain things, and put
the manuscript aside. Various events and projects distracted their
attention and the manuscript remained unsubmitted for over a decade.
In March 2018, Ruth died of lung cancer. In 2019, at the urging of 
guest editor Marius Tucsnak, the remaining author has made minimal
adjustments and has submitted the manuscript to MCRF. The references
remained as they were in 2007, except that papers with the status of
``submitted'' or ``accepted'' have been updated. 

\section{Some background on system nodes and well-posed linear
               systems} 

\ \ \ In this section, we recall a rather general class of
infinite-dimensional linear systems, called system nodes. System nodes
do not satisfy any well-posedness assumption, but nevertheless they
have well defined state trajectories and output functions
corresponding to smooth input functions and compatible initial states,
see Proposition \ref{NonWellPosedSystemProp}. We also recall a few
facts about well-posed linear systems.

In the semigroup approach to infinite-dimensional systems, we often
encounter an operator semigroup $\tline$ on a Hilbert space $X$ which
determines two additional Hilbert spaces, denoted by $X_1$ and
$X_{-1}$. This construction is now part of standard operator semigroup
theory, see for example Engel and Nagel \cite{EnNa}, Staffans
\cite[Section 3.6]{stafbook} or Weiss \cite[Section 3]{weiss1}, so
that we recall the main facts without proof:

\begin{proposition} \label{RiggedHilbertSpaceProp}
Let $X$ be a Hilbert space and let $A\colon\Dom{A}\rarrow X$ be the
infinitesimal generator of a strongly continuous semigroup
$\tline=(\tline_t)_{t\geq 0}$ on $X$. Take $\alpha\in\rho(A)$.

\noindent $(1)$ For each $x\in\Dom A$, define $\|x\|_1=\|(\alpha I-A)
x\|$. Then $\|\cdot\|_1$ is a norm on $X_1$ which makes $X_1$ into a
Hilbert space, and $A\in\Lscr(X_1,X)$. The operator $(\alpha I-A)
^{-1}$ maps $X$ isometrically onto $X_1$.

\noindent $(2)$ Let $X_{-1}$ be the completion of \m $X$ with respect
to the norm $\|x\|_{-1}=\break \|(\alpha I-A)^{-1}x\|$. Then $X$ is a
Hilbert space and $A$ has a unique extension to an operator
$A\in\Lscr(X, X_{-1})$. $(\alpha I-A)^{-1}$ maps $X_{-1}$
isometrically onto $X$.

\noindent $(3)$ The restrictions of \m $\tline_t$ to $X_1$ form a
strongly continuous semigroup on $X_1$. The generator of \m $\tline$
as a semigroup on $X_1$ is the restriction of $A$ to $\Dom{A^2}$.

\noindent $(4)$ The operators $\tline_t$ have a unique extension to
$X_{-1}$ which form a strongly continuous semigroup on $X_{-1}$. The
generator of this extended semigroup is $A$ extended to $X$ as in
point $(2)$.

\noindent $(5)$ The choice of $\alpha\in\rho(A)$ does not change the
spaces $X_1$ or $X_{-1}$, since different values of $\alpha\in\rho(A)$
lead to equivalent norms on $X_1$ and $X_{-1}$.
\end{proposition}

The notations $A$ and $\tline_t$ will be used also for the various
restrictions and extensions of these operators (as in the above
proposition).

Below we give what we think to be the simplest formulation of the
definition of a system node. Equivalent and related definitions can be
found, for example, in Opmeer \cite{Opm}, Malinen {\em et al}
\cite{MalStafWei}, and Staffans \cite{StafPas,stafbook}. We refer to
the same sources for the proofs of the results stated in this section.

\begin{definition} \label{SysNodeDef} {\rm
Let $U$, $X$ and $Y$ be Hilbert spaces. Let $A$ be the infinitesimal
generator of a strongly continuous semigroup $\tline$ on $X$, $B\in
\Lscr(U,X_{-1})$ and $C\in\Lscr(X_1,Y)$. Let $\GGG:\rho(A)\rarrow
\Lscr(U,Y)$ be such that
\BEQ{TransferFunctDiff}
   \GGG(s) - \GGG(\beta) \m=\m C[(sI-A)^{-1} - (\beta I-A)^{-1}]B \m,
\EEQ
for all $s,\beta\in\rho(A)$. Then $\Sigma=(A,B,C,\GGG)$ is a {\em
system node} on $(U,X,Y)$.

$U$ is the {\em input space} of $\Sigma$, $X$ is its {\em state
space}, $Y$ is its {\em output space}, $A$ is its {\em semigroup
generator}, $B$ is its {\em control operator}, $C$ is its {\em
observation operator} and $\GGG$ is its {\em transfer function}.
$(A,B,C)$ is the {\em generating triple} of \m $\Sigma$.}
\end{definition}

For a system node $\Sigma$ it is useful to introduce the space 
\BEQ{Vdef}
   V \m=\m \left\{ \bbm{x\\ v}\in X\times U\ |\
   Ax+Bv\in X \right\} \m,
\EEQ
which is a Hilbert space with the norm
$$ \left\|\bbm{x\\ v}\right\|_V^2 \m=\m \|x\|^2 +
   \|v\|^2 + \|Ax + Bv\|^2 \m.$$
We recall the space $\Dscr_0$ introduced in Section 1:
\BEQ{Do}
   \Dscr_0\m=\m \left\{ \bbm{x\\ v}\in X\times C^2([0,\infty),U)\ |\
   \bbm{x\\ v(0)}\in V \right\}.
\EEQ
      
The operator $\CD\in\Lscr(V,Y)$, called the {\em combined
observation/feedthrough operator} of $\Sigma$, is defined by
$$ \CD \bbm{x \\ v} \m=\m C \left[x - (\beta I-A)^{-1}B v \right]
   + \GGG(\beta) v \m.$$
This makes sense, because for $\sbm{x\\ v}\in V$ we have \m $x-(\beta
I-A)^{-1}Bv\in\Dscr(A)$. It is easy to verify (using \rfb
{TransferFunctDiff}) that the above definition of $\CD$ is independent
of the choice of $\beta\in\rho(A)$. We have $\GGG(s)=\CD\sbm{(sI-A)
^{-1}B\\ I}$. To make the connection with more familiar formulas from
finite-dimensional linear systems theory, note the following
particular case: if $B\in\Lscr(U,X)$ and $C\in\Lscr (X,Y)$, then the
limit $D=\lim_{\l\rarrow\infty}\GGG(\l)$ exists and we have $\CD=[C\
D]$, $\GGG(s)=C(sI-A)^{-1}B+D$.

Suppose that $\Sigma$ is a system node on $(U,X,Y)$ with generating
triple $(A,B,C)$ and transfer function $\GGG$. Choose $\beta\in
\rho(A)$ and introduce the Hilbert space \m $Z=\Dscr(A)+(\beta
I-A)^{-1}BU$, with the norm
$$ \|z\|^2_Z \m=\m \inf\{\|x\|^2_1+\|v\|^2_U \ | \ \
   x\in\Dscr(A), \ v\in U, \  z=x+(\beta I-A)^{-1}Bv\}.$$
(The space $Z$ is independent of the choice of $\beta$.) The system
node $\Sigma$ is called {\it compatible} if $C$ has an extension
$\overline{C}\in\Lscr(Z,Y)$.

Such an extension $\overline{C}$ is usually not unique. If $\Sigma$ 
is compatible, then for every extension $\overline{C}$ we can find a
unique $D\in\Lscr(U,Y)$ such that 
$$ \CD \bbm{x \\ v} \m=\m \overline{C}x+Dv 
   \FORALL \bbm{x \\ v}\in V \m.$$
Hence, in this case, $\GGG(s)=\overline{C}(sI-A)^{-1}B+D$ for every
$s\in\rho(A)$.

The following proposition from \cite[Sect.~2]{MalStafWei} shows that
at least for suitably smooth input fuctions and compatible initial 
states, $\Sigma$ has well defined state trajectories and output 
functions.

\begin{proposition} \label{NonWellPosedSystemProp}
Let $\Sigma$ be a system node on $(U,X,Y)$, let $A,B$ and $V$ be as in
\rfb{Vdef} and let $\CD$ be the combined observation/feedthrough
operator of $\Sigma$. For all $\sbm{z_0\\ u}\in \Dscr_0$ the equation
\BEQ{SysEqs}
   \dot z(t) \m=\m Az(t) + Bu(t) \m,\qquad \qquad z(0) = z_0,
\EEQ
has a unique (classical) solution satisfying $z\in C^1([0,\infty),X)
\cap C^2([0,\infty),X_{-1})$, $\sbm{z \\ u}\in C([0,\infty),V)$. The
corresponding output function $y$, defined by
$$y(t) \m=\m \CD \bbm{z(t)\\ u(t)} \m,$$
is in $C([0,\infty),Y)$.
\end{proposition}

It is easy to show that under the assumptions of Proposition
\ref{NonWellPosedSystemProp}, if $\ddot u(t)=O(e^{\o t})$ as $t\to
\infty$ for some $\o\in\rline$, then the Laplace transforms $\hat u$,
$\hat z$, and $\hat y$ satisfy
\begin{eqnarray} \label{Laplace1}
   \hat z(s) &=& (sI - A)^{-1}[z_0 + B\hat u(s)],\\
   \hat y(s) &=& C(sI - A)^{-1}z_0 + \GGG(s)\hat u(s),\label{Laplace2}
\end{eqnarray}
for all $s\in\cline$ for which $\Re s$ is larger than $\o$ and also
larger than the growth bound of $\tline$ (see, for example, \cite
[Lemma 4.7.11]{stafbook}). If $z_0\in X$ and the input signal $u$ is a
distribution with a Laplace transform defined on some right
half-plane, then \rfb{Laplace1} and \rfb{Laplace2} define the signals
$z$ and $y$ as distributions that have polynomially bounded Laplace
transforms on some right half-plane.

The {\it dual system node} $\Sigma^d$ has the generating triple
$(A^*,C^*,B^*)$ and transfer function $\GGG^d$ defined by
$\GGG^d(s)=\GGG(\bar{s})^*$.

We obtain the well-known class of \emph{well-posed linear systems} by
adding one more assumption to those in Definition \ref{SysNodeDef}.

\begin{definition} \label{WPLSDef} {\rm
Let $\Sigma$ be a system node on $(U,X,Y)$. We call $\Sigma$ {\em 
well-posed} if for some (hence for every) $t>0$ there exists an $M_t
\geq 0$ such that
\BEQ{WellPosednessBalance}
    \|z(t)\|^2 + \|y\|_{L^2([0,t],Y)}^2 \m\leq\m M_t
    \bigpar{\|z_0\|^2 + \|u\|_{L^2([0,t],U)}^2} \m,
\end{equation}
for all $z$, $y$, $z_0$, and $u$ satisfying the conditions of
Proposition \ref{NonWellPosedSystemProp}.}
\end{definition}

A well-posed system node is usually called a {\it well-posed linear
system}. Necessary and sufficient conditions for well-posedness were
given in Curtain and Weiss \cite{CurWei}. For alternative definitions,
background and examples we refer to Salamon \cite{sala}, Staffans
\cite{stafTAMS}, \cite{stafbook}, \cite{StWe}, Weiss \cite{weiss10},
\cite{weiss12}, Weiss and Rebarber \cite{WeRe} and Weiss, Staffans and
Tucsnak \cite{WST}. All well-posed systems are compatible, see
Staffans and Weiss \cite{StWe}.

If the system node $\Sigma$ is well-posed, then it defines a
family of bounded operators parameterized by $t\geq 0$,
\begin{eqnarray} \label{sigma}
   \Sigma_t \m=\m \left[\begin{array}{cc}\tline_t&\Psi_t\\
   \Phi_t & \fline_t \end{array} \right] \m
\end{eqnarray}
from $\left[\begin{array}{c}X\\ L^2([0,t];U)\end{array}\right]$ to
$\left[\begin{array}{c}X\\ L^2([0,t];Y)\end{array}\right]$, such that 
$$ \Sigma_t\left[\begin{array}{c}z_0\\ \PPP_t u\end{array}\right]
   \m=\m \left[\begin{array}{c}z(t)\\ \PPP_t y \end{array}\right]$$
for all $z_0,z(t),u$ and $y$ as in Proposition
\ref{NonWellPosedSystemProp}. Here, $\PPP_t u$ and $\PPP_t y$ are the
restrictions of $u$ and $y$ to $[0,t]$. In fact, well-posed linear
systems are usually defined via the operator family $(\Sigma_t)_{t
\geq 0}$ by imposing certain algebraic conditions.

The state trajectories of the well-posed $\Sigma$ are given by  
\BEQ{tlinePhi}
   z(t) \m=\m \tline_t\mm z_0 + \Phi_t\mm u \m,
\EEQ
where $u\in\LEloc{U}$ is the input function. The operators $\Phi_t$
are  called the {\it input maps} of \m $\Sigma$ and they are given by
\BEQ{Phi_t}
   \Phi_t u \m=\m \int_0^t \tline_{t-\sigma} B u(\sigma)\m d\sigma \m.
\EEQ
The output function of the well-posed \m $\Sigma$ is given by
\BEQ{output}
   y \m=\m \Psi\m z_0 + \fline \m u \,,
\EEQ
where $\Psi:X\rarrow\m\LEloc{Y}$ is called the (extended) {\it output
map} of $\Sigma$ and $\fline$ is a continuous linear operator from
$\LEloc{U}$ to $\LEloc{Y}$ called the (extended) {\it input-output
map} of $\Sigma$. $\Psi$ is given by
\BEQ{Psirep}
   (\Psi z_0)(t) \m=\m C\m\tline_t\m z_0 \FORALL z_0\in\Dscr(A) 
\EEQ
and for every $t\geq 0$ we have $\Psi_t=\PPP_t\Psi$. $\fline$ is given
by
$$(\fline u)(t) \m=\m \CD \bbm{\Phi_t u\\ u(t)} \m,$$
for all $u\in C^2([0,\infty),U)$ such that $u(0)=0$. For every $t\geq
0$ we have $\fline_t=\PPP_t\fline\PPP_t$. The transfer function of
every well-posed system is {\em proper}, meaning that it is uniformly
bounded on some right half-plane.
 
Let us now return to the more general case of a system node $\Sigma$
on $(U,X,Y)$, with input signal $u$ and output signal $y$. We are
interested in the feedback control law $u=Ky+v$, where $K\in\Lscr(Y,
U)$ and $v$ is the new input signal. This corresponds to Figure 1, but
with $-K$ in place of $\ka$. Elementary manipulations lead to
\BEQ{star}
   [I-\GGG(s) K]\hat y(s) \m=\m C(s I- A)^{-1}z_0 + \GGG(s)\hat v(s)
\EEQ 
and if the operators $I-K\GGG(s)$ have bounded inverses for all $s$ in
some right half-plane, then we can compute
\begin{eqnarray} \label{2star} 
   \hat z(s) &=&  (s I- A)^{-1}[z_0 + B
   K(I-\GGG(s) K)^{-1}C(sI-A)^{-1}z_0]\nonumber \\
   && +(sI-A)^{-1}BK(I-\GGG(s) K)^{-1}\GGG(s)\hat{v}(s) \m. 
\end{eqnarray}
In general, such a feedback need not result in a system node (even if
the inverse of $I-\GGG(s)K$ is uniformly bounded on some right
half-plane). In fact, a closed-loop semigroup generator may not exist.
In this paper we are interested in the special case where the feedback
results in a system node $\Sigma^K$, called the {\it closed-loop
system}. In this case we call $K$ an {\it admissible feedback
operator} for the system node $\Sigma$. A necessary condition for this
is the bounded invertibility of the operators $I-K\GGG(s)$ for all $s$
in some right half-plane. From \rfb{star} we see that the transfer
function of $\Sigma^K$ is
$$\GGG^K \m=\m \GGG(I-K\GGG)^{-1} \m=\m (I-\GGG K)^{-1}\GGG$$
and its observation operator $C^K$ is determined by
$$C^K(sI-A^K)^{-1} = (I-\GGG(s) K)^{-1}C(sI-A)^{-1} \m.$$
 From \rfb{2star} we see that the semigroup generator $A^K$ is
determined by
\BEQ{AK} 
       (sI-A^K)^{-1} \m=\m (s I- A)^{-1} + 
       (s I- A)^{-1}BK(I-\GGG(s) K)^{-1}C(sI-A)^{-1},
\EEQ
and the control operator $B^K$ is determined by
$$(sI-A^K)^{-1}B^K =  (s I- A)^{-1}B(I-K\GGG(s))^{-1} \m.$$
The last three formulas hold for all $s\in\rho(A)\cap\rho(A^K)$. The
latter is not empty since both $A$ and $A^K$ are semigroup generators.

In particular, if the closed-loop system $\Sigma^K$ is well-posed,
then the $\Lscr(U)$-valued function $(I-K\GGG)^{-1}$ is proper and we
call $K$ a {\it well-posed feedback operator} for $\Sigma$. Well-posed
feedback operators for well-posed systems were studied in
\cite{weiss12}, \cite{stafbook}.

\begin{proposition} \label{X}
Let $\Sigma$ be a system node on $(U,X,Y)$ and let $K\in
\Lscr(Y,U)$. Then $K$ is an admissible feedback operator for $\Sigma$
if and only if $K^*$ is an admissible feedback operator for the dual
system node $\Sigma^d$. Similarly, $K$ is a well-posed feedback
operator for $\Sigma$ if and only if $K^*$ is a well-posed feedback
operator for $\Sigma^d$.

If $K$ is an admissible feedback operator for $\Sigma$, resulting in
the closed-loop system $\Sigma^K$, then the feedback operator $K^*$
for $\Sigma^d$ leads to the closed-loop system $\Sigma^{d K^*}$, which
is the dual of \m $\Sigma^K$.
\end{proposition}

\section{Impedance passive and scattering passive systems} 

\ \ \ In this section, we give rigorous definitions for the concepts
of scattering passive and impedance passive system nodes following the
terminology in \cite{StafPas,StafCol,StafMTNS02}. We start with the
much simpler case of discrete-time systems.

\begin{definition} \label{Jacky_speaks} {\rm
Let $U,X,Y$ be Hilbert spaces. A {\em discrete-time system} $\Sigma_d$
on \\ $(U,X,Y)$ consists of four operators $A_d\in\Lscr(X)$, $B_d\in
\Lscr(U,X)$, $C_d\in\Lscr(X,Y)$ and $D_d\in\Lscr(U,Y)$, also called
the {\em generating operators} of \m $\Sigma_d$. If $u$ is an {\em
input signal} for the system (an arbitrary $U$-valued sequence defined
for $k\in\{0,1,2,\ldots\}$), then a {\em state trajectory} of the
system corresponding to this input is a solution $z$ of the difference
equation $z_{k+1}=A_dz_k+Bu_k$. This space trajectory is uniquely
determined if we specify the {\em initial state} $z_0\in X$. The
corresponding {\em output signal} $y$ is the $Y$-valued sequence
defined by $y_k=Cz_k+Du_k$.}
\end{definition}

While the Laplace-transform is a common tool for analysing
continuous-time systems, the analogous tool for discrete-time systems
is the $\Zscr$-transform: if $u$ is a sequence then its {\em
$\Zscr$-transform} is $\hat u(z)=\sum_{k=0}^\infty u_k z^{-k}$. If the
series converges for some $\zeta\in\cline$, then it converges for
every $z\in\cline$ with $|z|>|\zeta|$. In this case, we call $u$ {\em
$\Zscr$-transformable}. With the notation of the last definition, if
$u$ is $\Zscr$-transformable and $z_0=0$, then also $y$ is
$\Zscr$-transformable and we have $\hat y(z)=\GGG_d(z)\hat u(z)$,
for all $z\in\cline$ with $|z|$ sufficiently large, where
$$\GGG_d(z) \m=\m C_d (zI-A_d)^{-1} B_d + D_d \FORALL z\in\rho(A)
            \m.$$
The $\Lscr(U,Y)$-valued function $\GGG_d$ is called the {\em transfer
function} of \m $\Sigma_d$.

\begin{definition} {\rm
With the notation of the last definition, $\Sigma_d$ is {\em
scattering passive} if for all $z_0\in X$, all input signals $u$
and all $m\in\nline$,
\BEQ{scatdis}
   \|z_m\|^2-\|z_0\|^2 \m\leq\m \sum^{m-1}_{k=0} \|u_k\|^2 -
   \sum ^m_{k=0}\|y_k\|^2 \m.
\EEQ 

$\Sigma^d$ is {\em impedance passive} if for all $z_0\in X$, all
input signals $u$ and all $m\in\nline$,}
\BEQ{impdis}
   \|z_m\|^2-\|z_0\|^2 \m\leq\m 2\Re\sum^{m-1}_{k=0} \langle u_k,
   y_k \rangle \m.
\EEQ
\end{definition}

In fact, it enough to verify \rfb{scatdis} or \rfb{impdis} for $m=1$,
as it is easy to see. Hence, \rfb{scatdis} holds if and only if 
$\bbm{A_d&B_d\\ C_d&D_d}$ is a contraction. In this case, $A_d$ is a
contraction and $\GGG_d$ satisfies $\|\GGG_d(z)\|\leq 1$ for all
$z\in\cline$ with $|z|>1$.

It is readily verified that \rfb{impdis} is satisfied if and only if
the following holds
\BEQ{disimp}
   \bbm{A_d^*A_d&A_d^*B_d\\ B_d^*A_d&B_d^*B_d} \m\leq\m
   \bbm{I & C_d^* \\ C_d & D_d+D_d^*} \m.
\EEQ
In this case, $A_d$ is (again) a contraction and $\GGG_d$ satisfies
$\GGG_d(z)+\GGG_d(z)^*\geq 0$ for all $z\in\cline$ with $|z|>1$.

\begin{definition} {\rm
With the notation of Definition \ref{SysNodeDef}, the {\em (internal)
Cayley transform} of \m $\Sigma$ for the parameter $\alpha\in\rho(A)$
is the following operator in $\Lscr(X\times U,X\times Y)$:}
\BEQ{Cayley}
   \left[\begin{array}{cc} {\mathfrak A}(\alpha) & {\mathfrak B}
   (\alpha)\\ {\mathfrak C}(\alpha) & \GGG(\alpha) \end{array}\right]
   \m=\m \left[ \begin{array}{cc} (\overline{\alpha}I+A)(\alpha I-
   A)^{-1} & \sqrt{2\Re\alpha}\m (\alpha I-A)^{-1}B\\ \sqrt{2\Re
   \alpha} \m C(\alpha I-A)^{-1} & \GGG(\alpha)
   \end{array}\right].
\EEQ
\end{definition}

Note that the Cayley transform of $\Sigma$ defines a discrete-time
linear system $\Sigma_d$ with generating operators $A_d={\mathfrak A}
(\alpha),\m B_d={\mathfrak B}(\alpha),\m C_d={\mathfrak C}(\alpha),\m
D_d={\GGG}(\alpha)$. Some algebraic computation (using
\rfb{TransferFunctDiff}) shows that the transfer function of \m
$\Sigma_d$ is
$$ \GGG_d(z) \m=\m \GGG\left( \frac{\alpha z-\overline{\alpha}}{z+1}
   \right) \FORALL z\in\rho(A_d) \m=\m \left\{ \left.\frac{\overline
   {\alpha}+\l}{\alpha-\l} \ \right|\ \l\in\rho(A) \right\} \m.$$
The above formula for $\rho(A_d)$ refers to the usual case when $A$ is 
unbounded. If $A$ is bounded, then $\rho(A_d)$ contains in addition
the point $-1$.

We denote by $\cline_0$ the open right half-plane:
$$\cline_0 \m=\m \{s\in\cline\ | {\rm Re}\ s>0\} \m,$$
and we denote by $\Hscr^2(\cline_0)$ the usual Hardy space of analytic
functions on $\cline_0$, see Rudin \cite{Rudin}. If $U$ is a Hilbert
space, then the Hardy space $\Hscr^2(\cline_0;U)$ (containing
analytic $U$-valued functions) is defined similarly, see for example
Nagy and Foias \cite{NagFoi} or Rosenblum and Rovniak \cite{RoRo}. The
inner product on $\Hscr^2(\cline_0;U)$ is
$$ \langle v,w\rangle \m=\m \frac{1}{2\pi} \int_{-\infty}^\infty
   \langle v(\jw),w(\jw) \rangle \dd\m\o \m,$$
where the values $v(\jw),\m w(\jw)$ are obtained as nontangential
limits, for almost every $\o\in\rline$. The {\em Paley-Wiener theorem}
states that the Laplace transformation is a unitary operator from
$\LE U$ \m to \m $\Hscr^2(\cline_0;U)$.

Take $\alpha\in\cline_0$. The functions $\varphi_k:\cline_0\rarrow
\cline$ defined by 
$$ \varphi_k(s) \m=\m \frac{\sqrt{2\Re\alpha}}{\overline{\alpha}+s}
   \left(\frac{\alpha-s}{\overline{\alpha}+s}\right)^k \m,$$
where $k\in\{0,1,2,\ldots\}$, form an orthonormal basis in $\Hscr^2
(\cline_0)$. (Indeed, for $\alpha=1$ this is well-known, and for other
$\alpha$ it follows by scaling and vertical shifting.) This implies
that every $\hat u\in\Hscr^2(\cline_0;U)$ can be written in a unique
way as
\BEQ{Laguerre}
   \hat u(s) \m=\m \sum_{k=0}^\infty u_k \varphi_k(s) \m,
\EEQ
where $u_k\in U$ and \vspace{-3mm}
$$ \|u\|^2 \m=\m \sum_{k=0}^\infty \|u_k\|^2 \m.$$

\begin{remark} {\rm
One fascinating fact in infinite-dimensional systems theory is that
the internal Cayley transformed system $\Sigma^d$ inherits many
important properties of $\Sigma$, and it also works in the opposite
direction: many properties of \m $\Sigma$ follow from properties of
$\Sigma^d$. The reason for this is as follows:
 
Let \m $\Sigma$ be a system node. Assume that its initial state is
zero and
$$u \m\in\m C^2([0,\infty),U)\cap\LE U$$
is such that the corresponding output signal $y$ (defined as in
Proposition \ref{NonWellPosedSystemProp}) satisfies $y\in\LE Y$. Then
(according to \rfb{Laplace2}) we have $\hat y(s)=\GGG(s)\hat u(s)$,
for all $s$ in some right half-plane. Using the representation
\rfb{Laguerre} for $\hat u$, and a similar one for $\hat y$, we obtain
that for $\Re s$ sufficiently large,
$$ \frac{\sqrt{2\Re\alpha}}{\overline{\alpha}+s} \sum_{k=0}^\infty y_k
   \left(\frac{\alpha-s}{\overline{\alpha}+s}\right)^k \m=\m \GGG(s)
   \frac{\sqrt{2\Re\alpha}}{\overline{\alpha}+s} \sum_{k=0}^\infty u_k
   \left(\frac{\alpha-s}{\overline{\alpha}+s}\right)^k \m.$$
Denoting $z=\frac{\overline{\alpha}+s}{\alpha-s}$ (so that $s=\frac
{\alpha z-\overline{\alpha}}{z+1}$), we obtain
$$ \sum_{k=0}^\infty y_k z^{-k} \m=\m \GGG\left(\frac{\alpha z-
   \overline{\alpha}}{z+1}\right) \sum_{k=0}^\infty u_k z^{-k}
   \m=\m \GGG_d(z) \sum_{k=0}^\infty u_k z^{-k} \m.$$
Thus, regarding $(u_k)$ as the input signal to the Cayley transformed
system $\Sigma_d$, with initial state zero, the sequence $(y_k)$ will
be the corresponding output signal. In other words, the input-output
map of \m $\Sigma_d$ corresponds to the input-output map of \m
$\Sigma$ through the unitary transformation that maps $u$ into the
sequence $(u_k)$. Similar statements hold for the input-to-state
map and the initial state-to-output map.}
\end{remark}

In the following definition we use the notation
$\Dscr_0$ from \rfb{Do}.

\begin{definition} \label{passive} {\rm
Let $\Sigma$ be a system node with input space $U$, state space $X$
and output space $Y$. $\Sigma$ is called {\it scattering passive} if
for all $\bbm{x\\ u} \in \Dscr_0$ and every $\tau>0$,
$$ \|z(\tau)\|^2-\|z_0\|^2 \m\leq\m \int_0^\tau \|u(t)\|^2 \m \dd t
   - \int_0^\tau \|y(t)\|^2 \m \dd t \m.$$
Here the state trajectory $z$ and the output function $y(t)$ are
defined as in Proposition \ref{NonWellPosedSystemProp}. $\Sigma$ is
called {\it impedance passive} if $Y=U$ and for all $\tau>0$,
\rfb{statepassive} holds.}
\end{definition}

 From Definition \ref{WPLSDef} we see that scattering passive system
nodes are automatically well-posed and system stable. The semigroup of
any scattering or impedance passive system is a semigroup of
contractions. We refer to Staffans \cite{StafPas}, \cite{StafMTNS02},
\cite[Section 7]{StWe} and Weiss and Tucsnak \cite[Section 4]{WeTu}
for more details on scattering passive systems. The following theorem
is contained in \cite[Theorem 7.4]{StWe} and \cite[Theorems 3.3, 4.2]
{StafPas}. If $X$, $A$ and $X_1$ are as in Section 2, we denote by
$Z_{-1}$ the dual of $X_1$ with respect to the pivot space $X$. The
spaces $U$ and $Y$ are identified with their duals. 

\begin{theorem} \label{tpassive}
Let $\Sigma,\m U,\m X$ and $Y$ be as in the last definition. We denote
the generating triple of the system node \m $\Sigma$ by $(A,B,C)$ and
its transfer function by $\GGG$.

\begin{itemize}
\item[{\rm (a)}] $\Sigma$ is scattering passive if and only if for
some (hence, for every) $s\in\rho(A)$
$$ \left[\begin{array}{cc} A+A^* & (sI+A^*)(sI-A)^{-1}B\\ B^*
   (\overline sI-A^*)^{-1}(\overline sI+A) & B^*(\overline sI-A^*)
   ^{-1} 2(\Re s)(sI-A)^{-1} B \end{array} \right]$$
$$+ \left[ \begin{array}{cc} C^*C & C^*\GGG(s)\\ \GGG(s)^*C & \GGG
    (s)^*\GGG(s)\end{array} \right] \m\leq\m \left[ \begin{array}{cc}
    0 & 0\cr 0 & I \end{array} \right] \m.$$ 
\item[{\rm (b)}] $\Sigma$ is impedance passive if and only if \ $Y=U$
and for some (hence, for every) $s\in\rho(A)$
$$ \left[\begin{array}{cc} A+A^* & (sI+A^*)(sI-A)^{-1}B\\ B^*
   (\overline sI-A^*)^{-1}(\overline sI+A) & B^*(\overline sI-A^*)
   ^{-1} 2(\Re s)(sI-A)^{-1} B \end{array} \right]$$
$$ \leq \left[ \begin{array}{cc} 0 & C^*\\
   C & \GGG(s) + \GGG(s)^*\end{array} \right] \m.$$  
\item[{\rm (c)}] $\Sigma$ is impedance passive if and only $Y=U$ and
for some (hence, for every) $\alpha\in\cline_0$ the Cayley transform
defined in \rfb{Cayley} satisfies
$$ \left[\begin{array}{cc}{\mathfrak A}(\alpha)^*{\mathfrak A}(\alpha)
   &{\mathfrak A}(\alpha)^*{\mathfrak B}(\alpha) \\ {\mathfrak B}
   (\alpha)^* {\mathfrak A}(\alpha) & {\mathfrak B}(\alpha)^*
   {\mathfrak B}(\alpha) \end{array} \right] \m\leq\m \left[
   \begin{array}{cc} I & {\mathfrak C}(\alpha)^* \\ {\mathfrak C}
   (\alpha) & \GGG(\alpha) + \GGG(\alpha)^* \end{array}\right]$$
\end{itemize}

Note that all the above matrices in {\rm (a)} and {\rm (b)} are
self-adjoint operators in $\Lscr(X_1\times U,Z_{-1}\times U)$, so
that they determine symmetric quadratic forms on $X_1\times U$. The
inequalities in {\rm (a), (b)} are understood in the sense of these
quadratic forms.
\end{theorem}

\begin{remark} \label{Aug_2_2007} {\rm
By point (a) of the above theorem, the transfer function $\GGG$ of any
scattering passive system is defined on $\cline_0$ and it satisfies
$\|\GGG(s)\|\leq 1$ for all $s\in\cline_0$.

By point (b) of the above theorem the transfer function of any
impedance passive system is necessarily {\it positive}, i.e., 
$$\GGG(s)+\GGG(s)^*\geq 0 \ \ \mbox{~for all~}\  s\in \cline_0 \m.$$
Both converse statements are false, even in finite dimensions. This is
because impedance or scattering passivity depends on the realisation
of $\GGG$.}
\end{remark}

\begin{remark} \label{energypreserving} {\rm
If equality holds in \rfb{statepassive}, then the system is called
{\em impedance energy preserving} and the condition in part (b) of
Theorem \ref{tpassive} becomes an equality. In this case $A$ generates
an isometric semigroup.}
\end{remark}

We recall that if $E\in\Lscr(U)$ is self-adjoint, then $U$ has a
unique orthogonal decomposition into $E$-invariant subspaces
$U=U^-\oplus U^+$ such that $\langle Ev,v\rangle\leq 0$ for all $v\in
U^-$ and $\langle Ev,v\rangle>0$ for all non-zero $v\in U^+$. The
operator $E^+$ obtained by redefining $E$ to be zero on $U^-$ is
called {\it the positive part} of $E$. Noting that $\|E\| I\geq
E^+\geq E$, we deduce the following corollary.

\begin{corollary} \label{cpassive} 
Let $\Sigma$ be a system node with transfer function $\GGG$. Let
$\Sigma_E$ be the system node with the same generating operators and
the transfer function $\GGG+E$, where $E=E^*\in \Lscr(U)$. If
$\Sigma_E$ is impedance passive, then for all $c\geq ||E^+||$ the
system node $\Sigma_{cI}$ is impedance passive.
\end{corollary}

\begin{remark} \label{remABB*} {\rm 
For the special case that $A$ is skew-adjoint and $C=B^*$, we see that
the necessary and sufficient condition for impedance passivity reduces
to
\BEQ{namur}
   \GGG(s)+\GGG(s)^*-B^*(\bar{s}I-A)^{-1}(2\Re s)(sI-A)^{-1}B
   \m\geq\m 0 \m,
\EEQ
for some (hence, for every) $s\in\rho(A)$. In particular, if there
exists an $\o\in\rline$ such that $i\o\in\rho(A)$, then \rfb{namur}
reduces to $\GGG(i\o)+\GGG(i\o)^*\geq 0$.}
\end{remark}

\begin{remark} \label{boundedpassive} {\rm 
As in Staffans \cite[p.~294]{StafPas}, it is easily verified that the
necessary and sufficient condition in (b) for a well-posed system with
a bounded generating triple $A,B,C$ to be impedance passive reduces to
$$ \left[\begin{array}{cc}-A-A^* & C^*-B\\C-B^* & D+D^*\end{array}
   \right] \m\geq\m 0 \m,$$ 
where  $D=\lim_{\l\to+\infty}\GGG(\l)$. In the case that $D+D^*$ is
invertible, this is equivalent to $D+D^*>0$ and \m $A+A^*+(C^*-B)
(D+D^*)^{-1}(C-B^*)\leq 0$.}
\end{remark} 

We have the following simple alternative condition for impedance
passivity.

\begin{proposition} \label{ABCimp} 
Let $\Sigma$ be a system node with generating triple $(A,B,C)$ and
transfer function $\GGG$. Let $\Sigma_E$ be the system node with the
same generating operators and the transfer function $\GGG+E$, where
$E=E^*\in \Lscr(U)$. If $i\o\in \rho(A)$, then $\Sigma_E$ is impedance
passive if and only if
\BEQ{ABC}
   \left[\begin{array}{cc}- A_\o^{-1}-A_\o^{-*} &
   A^{-1}_\o B+A^{-*}_\o C^*\\  B^*A^{-*}_\o+C A^{-1}_\o  &
   2E+\GGG(i\o)+\GGG(i\o)^* \end{array}\right]
   \m\geq\m 0\m,
\EEQ
where $A_\o=A-i\o I$. 
\end{proposition}

{\it Proof.} \m  According to Theorem \ref{tpassive}(b) with $s=\jw$,
$\Sigma_E$ is impedance passive if and only if
$$ \left[\begin{array}{cc} A_\o+A_\o^* &
   -A_\o^*A_\o^{-1}B-C^*\\-B^*A_\o^{-*}A_\o-C&
   -2E-\GGG(i\o)-\GGG(i\o)^*\end{array}\right]
   \m\leq\m 0 \m.$$
Premultiplying this with the block operator diag$(A_\o^{-*},I)$ and
postmultiplying with diag$(A_\o^{-1},I)$ yields \rfb{ABC}. \kasten
 
\begin{remark} \label{remABC1} {\rm
Under the assumptions of Proposition \ref{ABCimp} define
$B_\o=A_\o^{-1} B$ and $C_\o=-CA_\o^{-1}$. Then \rfb{ABC} is the
necessary and sufficient condition for the {\it reciprocal system}
with bounded generating triple $(A_\o^{-1},B_\o,C_\o)$ and transfer
function
$$\GGG^r(s) \m=\m \GGG(i\o)+C_\o(sI-A_\o^{-1})^{-1}B_\o$$
to be impedance passive. If $2E+\GGG(i\o)+\GGG(i\o)^*$ is invertible,
then the necessary and sufficient condition \rfb{ABC} for the
impedance passivity of \m $\Sigma_E$ becomes
$$ A_\o^{-1}+A_\o^{-*}+(A_\o^{-1}B+A_\o^{-*}C^*)
   (2E+\GGG(i\o)+\GGG(i\o)^*)^{-1}(B^*A_\o^{-*}+CA_\o^{-1})\leq 0.$$
In particular, if $A$ is skew-adjoint, then $A_\o^*=-A_\o$ and the
above necessary and sufficient condition reduces to $B=C^*$ and}
\BEQ{nec}
   2E+\GGG(i\o)+\GGG(i\o)^*\geq 0.
\EEQ
\end{remark}

An easy consequence of Proposition \ref{ABCimp} is the following.

\begin{corollary} \label{corASS} 
Suppose that \m $\Sigma$ is a system node with generating triple
$(A,B,C)$ and transfer function $\GGG$, where $A$ generates a
contraction semigroup on $X$ and there exists an $\o\in\rline$ such
that \m $\jw\in\rho(A)$ and
\BEQ{ASS}
   B^*(i\o I+A^*)^{-1} \m=\m C(i\o I-A)^{-1} \m.
\EEQ
Then $\Sigma$ is impedance passive if and only if
$$\GGG(\jw)+\GGG(\jw)^* \m\geq\m 0 \m.$$
Hence, the minimal self-adjoint $E$ for which $\Sigma_E$ is impedance
passive is
$$E=-\half[\GGG(i\o)+ \GGG(i\o)^*].$$
\end{corollary}
\medskip


An interesting class of systems with colocated actuators and sensors
was investigated in Oostveen \cite{Oos}: $A$ generates a contraction
semigroup, $B$ is bounded and $C=B^*$. It is well-known that under
these assumptions $\GGG(s)=B^*(sI-A)^{-1}B$ is positive and below we
show that the system is also impedance passive.

\begin{proposition}\label{job} 
Let $\Sigma$ be well-posed linear system on $(U,X,U)$ with generating
operators $(A,B,B^*)$, where $B\in\Lscr(U,X)$ and $A$ generates a
contraction semigroup on $X$. Then $\Sigma$ is impedance passive.
\end{proposition}

{\it Proof.} \m Rearranging the terms in part (b) of Theorem \ref
{tpassive} we see that \m $\Sigma$ is impedance passive if and only
if for some $s\in\rho(A)$,
\BEQ{job1}
   \left[\begin{array}{cc} A+A^* &(A+A^*)(sI-A)^{-1}B \\ 
   B^*(\bar{s} I-A^*)^{-1}(A+A^*)  & B^*(\bar{s}I-A^*)^{-1}(A+A^*)
   (sI-A)^{-1}B \end{array}\right]
   \m\leq\m 0\m,
\EEQ
and this can be factored as 
\BEQ{matrixineq}
   \left[\begin{array}{cc} I &0 \\ 0 & B^*(\bar{s} I-A^*)^{-1}
   \end{array} \right] \left[\begin{array}{cc} A+A^* & A+A^*\\ A+A^*
   & A+A^* \end{array} \right] \left[\begin{array}{cc}
   I & 0 \\ 0 & (s I-A)^{-1}B \end{array} \right] \m\leq\m 0 \m.
\EEQ
Since $B$ is bounded, the left-hand side defines a symmetric quadratic form on $X_1\times U$ and it is clearly non-negative. \kasten

In Curtain and Weiss \cite{CuWeABB*} it was shown by means of
counter-examples that for $B$ unbounded the above proposition is
false. So it is interesting to ask if it is possible to find an
explicit expression for the minimal $E$ to ensure that $\Sigma_E$ is
impedance passive. In \cite{CuWeABB*} this was done for the special
class of well-posed systems satisfying the following two assumptions.

{\bf Assumption ESAD.} The operator $A$ is {\it essentially
skew-adjoint and dissipative}, which means that $\Dscr(A)=\Dscr(A^*)$
and there exists a $Q\in\Lscr(X)$ with $Q\geq 0$ such that
\BEQ{ESAD}
   Ax + A^*x \m=\m -Qx \FORALL x \in \Dscr(A) \m.
\EEQ

{\bf Assumption COL.} $Y=U$ and $C=B^*$.

\begin{remark} {\rm 
For a well-posed system satisfying the above two assumptions, in
\cite{CuWeABB*} (in Theorem 5.2 and its proof) it was shown that the
expression
\BEQ{Efirst}
   E \m=\m -\half[\GGG(s)+\GGG(s)^*] +\half B^* (\bar{s}I-A^*)^{-1}
   \left[2({\rm Re}s) I+Q \right](sI-A)^{-1}B \m,
\EEQ
is independent of $s\in\cline_0$ and this operator can be written
also in the form
$$ E \m=\m -\half \lim_{\l\rarrow\infty} \left[ \GGG(\l)^*
   +\GGG(-\l)\right] \m.$$
Under an additional smoothness assumption on $B$, it was shown in
Proposition 5.5 of \cite{CuWeABB*} that the above $E$ is the minimal
self-adjoint operator ensuring the positivity of the transfer function
$\GGG+E$.}
\end{remark}

Now we show that the expression \rfb{Efirst} is independent of $s$
also for a system node with colocated actuators and sensors and, with
this operator $E$, the system node $\Sigma_E$ is impedance
passive. Here we allow $s\in\rho(A)$.
 
\begin{theorem} \label{syspass} 
Let $\Sigma$ be a system node satisfying the assumptions {\bf ESAD}
and {\bf COL} and let $E$ be defined by \rfb{Efirst} where $s\in\rho
(A)$. Then $E$ is independent of $s\in\rho(A)$ and the system node
$\Sigma_E$ with generating triple $(A,B,B^*)$ and transfer function
$\GGG_E=\GGG+E$ is impedance passive. Moreover, $E$ is the smallest
self-adjoint operator in $\Lscr(U)$ which has this passivity property.
\end{theorem}

{\it Proof.} \m First we note that the proof of Theorem 5.2 in \cite
{CuWeABB*} that the right-hand side of \rfb{Efirst} is a constant also
holds for a system node and for all $s\in\rho(A)$. For every $F=F^*\in
\Lscr(U)$, let us denote by $\Sigma_F$ the system node with generating
triple $(A,B,B^*)$ and transfer function $\GGG_F=\GGG+F$. For
simplicity we prove the result for a real $s=\l$. According to Theorem
\ref{tpassive}, $\Sigma_F$ is impedance passive if and only if for
some real $\l>0$,
\BEQ{Olof}
   \left[\begin{array}{cc} 2 \l(\l I-A^*)^{-1} Q (\l I-A)^{-1} &
   \sqrt{2\l}(\l I-A^*)^{-1}Q(\l I-A)^{-1}B\\ \sqrt{2\l} B^*
   (\l I-A^*)^{-1}Q(\l I-A)^{-1}  & {\bf P}_F(\l)\end{array}\right]
   \m\geq\m 0\m,
\EEQ
where ${\bf P}_F(\l) \m=\m \GGG_F(\l)^*+\GGG_F(\l)-B^*(\l I-A^*)^{-1}
2\l(\l I-A)^{-1}B$. If \rfb{Olof} holds for some $\l>0$, then it holds
for all $\l>0$. Using \rfb{Efirst} we find that
\BEQ{PF}
   {\bf P}_F(\l) \m=\m 2(F-E) + B^*(\l I-A^*)^{-1}Q(\l I-A)^{-1}B \m.
\EEQ
We see from \rfb{Olof} that $\Sigma_F$ is impedance passive if and
only if for some $\l>0$,
\begin{eqnarray} \label{matrixineq_bis}
   \left[\begin{array}{cc} \sqrt{2\l}(\l I-A^*)^{-1} &0 \\ 
   0 & B^*(\l I-A^*)^{-1} \end{array} \right]
   \left[\begin{array}{cc} Q & Q \\ Q & Q \end{array} \right]
   \nonumber \\ \left[\begin{array}{cc}\sqrt{2\l}(\l I-A)^{-1} & 0 \\
   0 & (\l I-A)^{-1}B \end{array} \right] \geq 
   \left[\begin{array}{cc} 0 & 0 \\ 0&2(E-F) \end{array} \right]\m.
\end{eqnarray}
Again, if this holds for some $\l>0$, then it holds for all $\l>0$.
It is clear that the left-hand expresssion is nonnegative and so we
see that \rfb{matrixineq_bis} holds for $F=E$, and that this is the
smallest possible value for which $\Sigma_F$ is impedance passive.
\kasten

\begin{remark} {\rm
If $Q=0$ ($A$ is skew-adjoint) and $F=E$, then clearly we have
equality in \rfb{Olof}. According to \cite[Theorem 4.6]{StafPas}, it
follows that if $Q=0$, then $\Sigma_E$ is impedance energy preserving
(see Remark \ref{energypreserving}).}
\end{remark}

\begin{remark}\label{boundedB} {\rm 
In the case that $B$ is bounded and $\GGG(s)=D+B^*(sI-A)^{-1}B$, then
taking limits as $\Re s\to \infty$ shows that $-2E=D+D^*$.}
\end{remark}

Another interesting class of systems is those with self-adjoint $A$.

\begin{lemma} \label{selfadjoint} 
Suppose that $\Sigma$ is a system node and that $A\leq 0$ (i.e., $A$
is self-adjoint and $\langle Ax,x\rangle\leq 0$ for all
$x\in\Dscr(A)$).
\begin{enumerate}
\item If $0\in\rho(A)$ and $C=B^*$, then the system node $\Sigma_E$
will be impedance passive if and only if
$$ 2E+\GGG(0)+\GGG(0)^*\geq 0 \m.$$
\item $\Sigma_E$ is impedance passive if for some $\alpha>0$
$$C^* \m=\m (\alpha I+A)(\alpha I-A)^{-1}B \m,$$
$$ 2E+\GGG(\alpha)+\GGG(\alpha)^*-2\alpha B^*(\alpha I-A)^{-2}B
   \m\geq\m 0 \m.$$
\item If $C^*=B$ and $(\alpha I-A)^{-\half}B$ is a bounded operator
for some positive $\alpha$, then there exists a minimal self-adjoint
$E$ for which $\Sigma_E$ is passive. It is given by
\BEQ{Esecond}
   E \m=\m-\half [\GGG(s)+\GGG(s)^*]+B^* (\bar{s}I-A)^{-1} [(\Re s)
   I-A](sI-A)^{-1}B \m, \EEQ which holds for all $s\in \cline_0$.
\end{enumerate}
\end{lemma}

{\it Proof.} \m We leave this to the reader as it is similar to
previous proofs.

\section{Impedance passive systems with feedback} 

\ \ \ The next result describes a feedback transformation from
impedance passive nodes to scattering passive well-posed systems. It
is closely related to \cite[Theorem 5.2]{StafPas} and to \cite
[Section 2]{weiss63}.
\begin{proposition} \label{scat} 
Let $\Sigma^p$ be an impedance passive system node on $(U,X,U)$ with
generating triple $(A,B,C)$ and transfer function $\GGG^p$ and let
$k>0$. Suppose that $u^p$ is a $C^2$ input function and $z_0$ an
initial state satisfying $Az_0+Bu^p(0)\in X$, and $z,\,y^p$ are the
resulting state trajectory and output function, as in Proposition
{\rm \ref{NonWellPosedSystemProp}}.
Then there is a unique scattering passive system node $\Sigma^s$ with
the following property: if we use the same initial state $z_0$ and the
new input function \ $u^s=\sqrt {\frac{k}{2}}(\frac{u^p}{k}+y^p)$ for
the system node $\Sigma^s$, then the state trajectory of \m $\Sigma^s$
is again $z$ and the output function of \m $\Sigma^s$ is given by \
$y^s=\sqrt{\frac{k}{2}}(\frac {u^p}{k}-y^p)$. We have \m
$(I+k\GGG^p)^{-1}\in H^\infty$ and the transfer function of \m
$\Sigma^s$ is \vspace{-1mm}
$$\GGG^s \m=\m (I-k\GGG^p)(I+k\GGG^p)^{-1} \m.$$
\end{proposition}

$\Sigma^s$ can be thought of as a closed-loop system node obtained
from $\Sigma^p$ as shown in Figure 2. In \cite{StafPas}, $k$ is chosen
to be 1, while in \cite{weiss63}, $\frac 1 r$ appears in place of
$k$. Staffans in \cite{StafPas} calls the transformation from $\Sigma
^p$ to $\Sigma^s$ the {\em diagonal transformation} (or the external
Cayley transformation) and he traces it back to the work of M. Livsic
\cite{LivsBook}.

{\it Proof.} \m Noting that $\GGG^p$ is positive, according to
Proposition 2.1 in \cite{CuWeABB*} (with $c=0$), it follows that
$(I+k\GGG^p)^{-1}\in H^\infty$. Thus the closed-loop transfer function
$$\GGG^s \m=\m (I-k\GGG^p)(I+k\GGG^p)^{-1} \m=\m 2(I+k\GGG^p)^{-1}-I$$
is in $H^\infty$. Hence, inputs $u^s\in\LE U$ produce outputs $y^s\in
\LE U$. Moreover, it is easy to check that the scattering passivity 
property holds, since
\BEQ{uy}
   \|u^s(t)\|^2 - \|y^s(t)\|^2 \m=\m 2\Re \langle y^p(t),u^p(t)
   \rangle \m.
\EEQ
It remains to show that $\Sigma^s$ is a system node. To do this we use
the internal Cayley transform introduced in Section 3, in particular
in \rfb{Cayley}.  First we note that from Theorem \ref{tpassive} the
Cayley-transformed system node produces the discrete-time linear
system with generating operators $A_d,B_d,C_d,D_d$ that satisfy
\rfb{disimp}.  So the discrete-time system is impedance passive. As
noted above, $I+k\GGG(\alpha)^p=I+kD_d$ is invertible for all $k>0$
and $\alpha \in \cline^+$. Hence the closed-loop discrete-time system
obtained via the transformation is well-defined and its generating
operators are given by $A_d^s,B_d^s,C_d^s,D_d^s$ are given by
$$ A^s_d \m=\m A_d-B_d(I+kD_d)^{-1}kC_d, \quad B^s_d \m=\m \sqrt{2k}
   B_d (I+kD_d)^{-1} \m,$$
$$ C^s_d \m=\m -\sqrt{2k}(I+kD_d)^{-1}C_d, \quad D_d^s \m=\m (I+k
   D_d)^{-1}(I-kD_d)^{-1} \m.$$

Using \rfb{uy} it is easy to see that it is scattering passive. Hence
$\bbm{A_d^s&B_d^s\\ C_d^s&D_d^s}$ is a contraction. In particular,
this implies that $A^s_d$ is a contraction. So its inverse Cayley
transform will define the generator of a contraction semigroup if and
only if -1 is not an eigenvalue of $A^s_d$. Using a contradiction
argument we show that -1 cannot be an eigenvalue of $A^s_d$. Suppose,
on the contrary, that $A^s_dx=-x$ for some non-zero $x\in X$. Then
from the contraction property we obtain
$$ \|x\|^2+\|C^s_dx\|^2=\|A^s_dx\|^2+\|C^s_dx\|^2 \m\leq\m \|x\|^2.$$
Hence $0=C^s_dx=-\sqrt{2k}(I+kD_d)^{-1}C_dx$ and $C_dx=0$. But 
$$-x \m=\m A^s_dx \m=\m A_dx-B_dk(I+kD_d)^{-1}C_dx \m=\m A_dx$$
implies that -1 is an eigenvalue of $A_d$. Since we have assumed that
$A$ generates a contraction semigroup, we have arrived at a
contradiction. The inverse Cayley transform of the closed-loop
discrete-time system is $\Sigma^s$ and it has a (contraction)
semigroup generator. Thus $\Sigma^s$ is a scattering passive system 
node. \kasten

\begin{figure}[htbp] 
\setlength{\unitlength}{0.2mm}
\begin{center}
\begin{picture}(670,130)
\put(  6, 55){\vector(1,0){54}}                
\put( 20, 63){\makebox(0,0)[lb]{{$u^s$}}}      
\put( 60, 25){\framebox(50,60){$\sqrt{\frac 2 k}$}} 
\put(110, 55){\vector(1,0){54}}                
\put(171, 55){\circle{14}}                     
\put(171,-12){\vector(0,1){60}}                
\put(146, 40){\makebox(0,0)[lb]{{$_+$}}}       
\put(175, 32){\makebox(0,0)[lb]{{$_-$}}}       
\put(178, 55){\vector(1,0){72}}                
\put(186, 63){\makebox(0,0)[lb]{{$\frac{u^p}{k}$}}}
\put(250, 25){\framebox(50,60){$k$}}           
\put(300, 55){\vector(1,0){54}}                
\put(320, 63){\makebox(0,0)[lb]{{$u^p$}}}      
\put(354, 25){\framebox(50,60){$\Sigma^p$}}    
\put(404, 55){\vector(1,0){62}}                
\put(424, 63){\makebox(0,0)[lb]{{$y^p$}}}      
\put(473, 55){\circle{14}}                     
\put(473,122){\vector(0,-1){60}}               
\put(458, 69){\makebox(0,0)[lb]{{$_+$}}}       
\put(448, 43){\makebox(0,0)[lb]{{$_-$}}}       
\put(480, 55){\vector(1,0){44}}                
\put(524, 25){\framebox(50,60){$\sqrt{\frac k 2}$}} 
\put(574, 55){\vector(1,0){54}}                
\put(596, 63){\makebox(0,0)[lb]{{$y^s$}}}      
\put(220, 55){\line(0,1){67}}                  
\put(220,122){\line(1,0){253}}                 
\put(436, 55){\line(0,-1){67}}                 
\put(171,-12){\line(1,0){265}}                 
\end{picture}
\end{center}
\caption{The scattering passive system $\Sigma^s$ with input $u^s$ and
output $y^s$, obtained from the impedance passive system node
$\Sigma^p$ via the diagonal transformation, as in Proposition 
\ref{scat}.
\hfill \m \hskip 90mm \m \hskip 37mm --------------------------}
\end{figure}
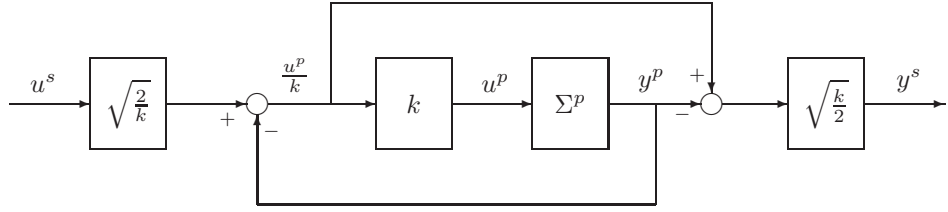

The relation between scattering passive and impedance passive system
nodes in Proposition \ref{scat} allows us to prove the following
stronger version of Theorem \ref{inpoutstab}.

\begin{theorem} \label{diss} 
Let $\Sigma$ be a system node on $(U,X,U)$ with generating triple
$(A,B,C)$ and transfer function $\GGG$. For any $E\in\Lscr(U)$ we
denote by $\Sigma_E$ the system node with the same generating triple
and with the transfer function $\GGG+E$. Assume that $E^*=E$ is such
that $\Sigma_E$ is impedance passive. Denote
$$c \m=\m \|E^+\| \m,\ \ \ \ \ka_0 \m=\m \frac{1}{c} \m,$$
where $E^+$ is the positive part of $E$ (see Section {\rm 3}).

Then for every $\ka\in (0,\ka_0)$, $K=-\ka I$ is a well-posed feedback
operator for \m $\Sigma$ and the corresponding closed-loop system \m
$\Sigma^\ka$ is system stable. Moreover, the semigroup of \m
$\Sigma^\ka$ is a contraction semigroup on $X$.
\end{theorem}

{\it Proof.} \m The generating triple of $\Sigma^\ka$ is complicated
to express directly in terms of $A,B,C,\GGG$ and $\ka$ (see for
example \cite{stafbook}). Instead, we use an indirect approach
involving several transformations from one system to another. Recall
that the system node $\Sigma$ has the transfer function $\GGG$. Let \m
$\Sigma^p$ be the impedance passive system node with transfer function
$\GGG^p=\GGG+cI$ from Corollary \ref{cpassive}. We apply Proposition
\ref{scat} to $\Sigma^p$, which implies that for every $k>0$, the
closed-loop system node $\Sigma^s$ with the transfer function
$\GGG^s=(I-k\GGG^p)(I+k\GGG^p)^{-1}$ is scattering passive and hence
well-posed. We have
\BEQ{ka} 
   I+k\GGG^p \m=\m (1+kc)I+k\GGG \m=\m (1+kc) \left( I+\frac{k}{1+kc}
   \m \GGG \right) \m.
\EEQ
Denote $\ka=\frac{k}{1+kc}$, so that $\ka\in(0,\ka_0)$. We know from
Proposition \ref{scat} that $(I+k\GGG^p)^{-1}\in H^\infty$, and this
together with \rfb{ka} shows that $(I+\ka\GGG)^{-1}\in H^\infty$. So
$\GGG^\ka=\GGG(I+\ka\GGG)^{-1}\in H^\infty$ and the input-output
connection in Figure 1 is well-posed. It remains to show that $-\ka I$
is a well-posed feedback operator for $\Sigma$ with the closed-loop
system node $\Sigma^\ka$. We do this by clarifying its relationship
with the three system nodes $\Sigma$, $\Sigma^p$ and $\Sigma^s$.  They
are all obtained from the same basic system node $\Sigma$ by defining
new input and output signals via linear transformations applied to the
original input signal $u$ and the original output signal $y$ of \m
$\Sigma$ (this does not affect the state trajectories). We start with
the connection with $\Sigma^\ka$: we denote its input by $v$ (as in
Figure 1). The output of $\Sigma^\ka$ is $y$, the same as for
$\Sigma$. Then $u=v-\ka y$ (see Figure 1), or in matrix form:
$$ \left[ \begin{array}{c} u\\ y\end{array} \right] \m=\m
   \left[ \begin{array}{cc} I & -\ka I\\ 0 & I\end{array} \right]
   \left[ \begin{array}{c} v\\ y\end{array} \right] \m.$$
The input and output signals of \m $\Sigma^p$ are $u^p=u$ and $y^p=y
+cu$. The input and output signals of \m $\Sigma^s$, denoted $u^s$ and
$y^s$, have been defined in terms of $u^p$ and $y^p$ in Proposition
\ref{scat} (see Figure 2). Writing these formulas in matrix form, we
have
$$ \left[ \begin{array}{c} u^p\\ y^p\end{array} \right] \m=\m
   \left[ \begin{array}{cc} I & 0\\ cI & I\end{array} \right]
   \left[ \begin{array}{c} u\\ y\end{array} \right] \m,\qquad
   \left[ \begin{array}{c} u^s\\ y^s\end{array} \right] \m=\m
          \sqrt{\frac k 2}
   \left[ \begin{array}{cr} \frac{1}{k} I & I\\ \frac{1}{k}I & -I
          \end{array} \right]
   \left[ \begin{array}{c} u^p\\ y^p \end{array} \right] \m.$$
To obtain the relation between the signals of $\Sigma^\ka$ and
$\Sigma^s$, we have to multiply the three $2\times 2$ matrices
appearing above, which yields
$$ \left[ \begin{array}{c} u^s\\ y^s\end{array} \right] \m=\m
          \sqrt{\frac k 2}
   \left[ \begin{array}{cr} (\frac{1}{k}+c)I & [1-(\frac{1}{k}+c)
          \ka] I \\ (\frac{1}{k}-c)I & -[1+(\frac{1}{k}-c)\ka] I
          \end{array} \right]
   \left[ \begin{array}{c} v\\ y \end{array} \right] \m=\m
   \left[ \begin{array}{cc} \frac{1}{\alpha}I & 0\\ \beta I &
          -\alpha I \end{array} \right]
   \left[ \begin{array}{c} v\\ y \end{array} \right] \m,$$
where \vspace{-1mm}
$$ \alpha \m=\m \sqrt{ 2\ka(1-\ka c)} \m,\qquad
   \beta  \m=\m \frac{1-2\ka c}{\alpha} \m,$$
so that $\alpha>0$. Inverting the last $2\times 2$ matrix, we have
\BEQ{Sigma_ka_from_d}
   \left[ \begin{array}{c} v\\ y \end{array} \right] \m=\m
   \left[ \begin{array}{cc} \alpha I & 0\\ \beta I
          & -\frac{1}{\alpha} I \end{array} \right]
   \left[ \begin{array}{c} u^s\\ y^s\end{array} \right] \m.
\EEQ
This relation is illustrated in Figure 3. \vspace{-2mm}

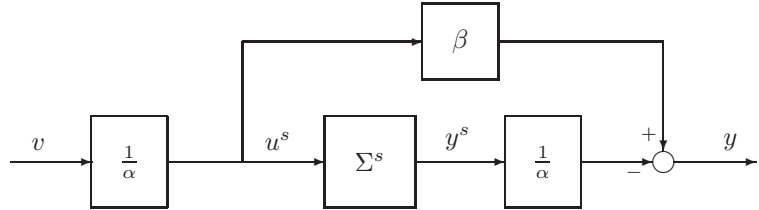
\begin{figure}[htbp]
\setlength{\unitlength}{0.2mm}
\begin{center}
\begin{picture}(700,140)
\put( 76, 25){\vector(1,0){54}}                
\put( 90, 33){\makebox(0,0)[lb]{{$v$}}}        
\put(130, -5){\framebox(50,60){$\frac 1 \alpha$}} 
\put(180, 25){\vector(1,0){105}}               
\put(350, 80){\framebox(50,50){$\beta$}}       
\put(245, 33){\makebox(0,0)[lb]{{$u^s$}}}      
\put(285, -5){\framebox(60,60){$\Sigma^s$}}    
\put(365, 33){\makebox(0,0)[lb]{{$y^s$}}}      
\put(495, 39){\makebox(0,0)[lb]{{$_+$}}}       
\put(485, 13){\makebox(0,0)[lb]{{$_-$}}}       
\put(345, 25){\vector(1,0){60}}                
\put(405, -5){\framebox(50,60){$\frac 1 \alpha$}} 
\put(455, 25){\vector(1,0){48}}                
\put(510, 25){\circle{14}}                     
\put(230, 25){\line(0,1){80}}                  
\put(230,105){\vector(1,0){120}}               
\put(400,105){\line(1,0){110}}                 
\put(510,105){\vector(0,-1){73}}               
\put(517, 25){\vector(1,0){55}}                
\put(550, 33){\makebox(0,0)[lb]{{$y$}}}        
\end{picture}
\end{center}
\caption{The feedback system $\Sigma^\ka$ from Figure 1, with input
$v$ and output $y$, as obtained from the scattering passive system
$\Sigma^s$ with input $u^s$ and output $y^s$. Note that there is no
feedback loop involved in this transformation.
\hfill \m \hskip 90mm \m \hskip 37mm --------------------------}
\end{figure}

We denote the generating triple of \m $\Sigma^s$ by $(A^s,B^s,C^s)$
and its transfer function by $\GGG^s$. Note from \rfb{Sigma_ka_from_d}
and Figure 3 that there is no feedback loop involved in the transform%
ation between the inputs and outputs of the scattering passive system
$\Sigma^s$ and of $\Sigma^\ka$. So the latter is in fact a well-posed
system and from \rfb{Sigma_ka_from_d} we deduce that its generating
triple $(A^\ka,B^\ka,C^\ka)$ and $\GGG^\ka$ are given by
\begin{eqnarray} \label{dk}
   A^\ka \m=\m A^s \m,\ \ \ B^\ka \m=\m \frac{1}{\alpha} B^s \m,\ \ \
   C^\ka \m=\m -\frac{1}{\alpha} C^s \m,\ \ \ \GGG^\ka \m=\m
   \frac{\beta}{\alpha} I - \frac{1}{\alpha^2} \GGG^s \m.
\end{eqnarray}
Since the scattering passive linear system $\Sigma^s$ is system
stable, the above relations imply that $\Sigma^\ka$ is also
system stable. Similarly, since $A^s$ generates a contraction
semigroup, so does $A^\ka$. \kasten

\section{Weak and strong stabilization} 

\ \ \ In this section, we continue our analysis of the stabilization
of impedance passive system nodes. We need the following result about
the spectrum of the closed-loop generator $A^K$ under admissible
feedback from Curtain and Jacob \cite[Theorem 6.2]{CurJac}, a
generalization of Lemma 4.4 in Salamon \cite{sala}, see also 
Weiss and Xu \cite{WeXu} for related results.

\begin{theorem} \label{pertnode} 
Let $\Sigma$ be a system node with generating triple $(A,B,C)$ and
transfer function $\GGG$. Suppose that $K\in\Lscr(Y,U)$ is an
admissible feedback operator for $\Sigma$ producing the closed-loop
system node $\Sigma^K$. Denote its generating triple by
$(A^K,B^K,C^K)$ and its transfer function $\GGG(I-K\GGG)^{-1}$ by
$\GGG^K$.

Then for $\l\in\rho(A)$ we have that $\l\in\rho(A^K)$ if and only if
$I-K\GGG(\l)$ is invertible.
\end{theorem}

Actually, in \cite{CurJac} this is proved for operator nodes, a
slightly more general concept than system nodes. The idea of the proof
is to apply the Cayley transform to $\Sigma$ and $\Sigma^K$ and then
to verify the corresponding discrete-time result.

For a system node $\Sigma$ on $(U,X,Y)$ with generating triple
$(A,B,C)$ we define the {\em unobservable space}
$$ \Nscr \m=\m \{x\in X\ |\ C(sI-A)^{-1}x=0 \mbox{~~for~~} \Re s >
   \o\} \m,$$ 
where $\o$ is some real number larger than the growth bound of the
semigroup generated by $A$. Clearly, the choice of $\o$ is
unimportant. Note that the approximate observability in infinite time
of \m $\Sigma$ is equivalent to $\Nscr=\{0\}$. The corresponding space
for the dual system node is
$$ \Nscr^d \m=\m \{x\in X\ |\ B^*(sI-A^*)^{-1}x=0 \mbox{~~for~~}
   \Re s >\o\} \m.$$ 
Clearly the approximate controllability in infinite time of \m
$\Sigma$ is equivalent to $\Nscr^d=\{0\}$. The main result of this
section is a strengthened version of Theorem \ref{approx_obs}:

\begin{theorem} \label{stability} 
Suppose that $\Sigma$ is a system node on $(U,X,U)$ with generating
triple $(A,B,C)$ and transfer function $\GGG$. For any $E\in\Lscr(U)$
we denote by $\Sigma_E$ the system node with the same generating
triple as $\Sigma$, but with the transfer function $\GGG+E$. Suppose
that there exists $E=E^*\in\Lscr(U)$ such that \m $\Sigma_E$ is
impedance passive, let $c=\|E^+\|$, where $E^+$ is the positive part
of $E$ and denote $\ka_0=\frac{1}{c}$.

Then for every $\ka\in(0,\ka_0)$ the operator $K=-\ka I$ is a
well-posed feedback operator for \m $\Sigma$. If either of the
following conditions holds,
\BEQ{Cweak}
   \left\{ x\in\Nscr \ |\ \|\tline_t x\|=\|x\|=\|\tline^*_t x\|
   \FORALL t>0 \,\right\} \m=\m \{0\} \m,
\EEQ
\BEQ{Bweak}
   \left\{ x\in\Nscr^d \ | \ \|\tline_t x\|=\|x\|=\|\tline^*_t x\|
   \FORALL t>0 \,\right\} \m=\m \{0\} \m,
\EEQ
then the semigroup $\tline^\ka$ of the closed-loop system \m 
$\Sigma^\ka$ is weakly stable.

Moreover, if \rfb{Cweak} or \rfb{Bweak} hold and $\sigma(A)\cap
i\rline$ is countable, then both semigroups $\tline^\ka$ and 
${\tline^\ka}^*$ are strongly stable.
\end{theorem}

{\it Proof.} \m (a) To prove weak stability, recall from Theorem \ref
{diss} that $\Sigma^\ka$ is a well-posed linear system and $\tline
^\ka$ is a semigroup of contractions. In Sz\"okefalvi-Nagy and Foias
\cite{NagFoi} or Davies \cite[Corollary 6.22]{Davi} (see also 
Benchimol \cite{Ben}) it is shown that $\tline^\ka$ is weakly stable
if $X^u=\{0\}$, where
$$ X^u \m=\m \{ x\in X \,|\ \|\tline^\ka_t x\|=\|x\|=\|\tline^{\ka*}_t
  x\| \FORALL t>0 \,\} \m,$$
i.e., $\tline^\ka$ is completely non-unitary. We show that \rfb{Cweak}
implies $X^u=\{0\}$. Theorem \ref{tpassive} applies to the scattering
passive system node $\Sigma^s$ introduced in the proof of Theorem
\ref{diss} (shown in Figure 2). Noting the close relationships
betweeen the generating operators of \m $\Sigma^\kappa$ and $\Sigma^s$
from \rfb{dk}, we obtain
$$ \langle A^\ka z,z\rangle + \langle {A^\ka}^*z,z \rangle +
   \alpha^2 \langle {C^\ka}^* C^\ka z,z\rangle \m\leq\m 0
   \FORALL z\in\Dscr(A^\ka) \m,$$
where $\alpha^2=2\ka(1-\ka c)>0$ and $\langle\cdot,\cdot\rangle$ is
the duality pairing between $\Dscr(A^\ka)$ and its dual with respect
to the pivot space $X$. We take $z=\tline^\ka_t x$ to obtain
$$ \frac{\dd}{\dd t} \|\tline^\ka_t \m x\|^2 +\alpha^2 \|C^\kappa\m
   \tline^\ka_t \m x\|^2 \m \leq \m 0 \FORALL x\in\Dscr(A^\ka) \m.$$
Integrating from 0 to $t$ gives \m $\|\tline^\ka_t x\|^2+\alpha^2\int
^t_0\|C^\kappa\tline^\ka_\sigma x\|^2\,\dd\sigma\leq\|x\|^2$, and this
can be extended by continuity to all of $X$:
\BEQ{weak1}
   \|\tline^\ka_t x\|^2+ \alpha^2\int^t_0 \| (\Psi^\kappa x)(\sigma)\|
   ^2 \,\dd \sigma \m\leq\m \|x\|^2 \FORALL x\in X \m,\ t\geq 0 \m.
\EEQ

Suppose now that $x\in X^u$. Substituting this $x$ into \rfb{weak1}
implies that $\Psi^\ka x=0$. Consider $\Sigma$ with the initial state
$x\in X^u$ and input $u=0$. Denote the corresponding output signal of
$\Sigma$ by $y$, which is defined by its Laplace transform
$$\hat{y}(s) \m=\m C(sI-A)^{-1}x \FORALL s\in\cline_0$$
(see Section 2). The initial state of $\Sigma^\ka$ is the same $x$ and
its input signal is $v=\ka y$ (see Figure 1). Its output signal
satisfies \m $y=\Psi^\ka x+\ka\fline^\ka y=\ka \fline^\ka y$. Taking
Laplace transforms we obtain $(I-\ka \GGG^\ka)\hat{y}=0$. Since
$I-\ka\GGG^\ka=(I+\ka \GGG)^{-1}$, it follows that $\hat{y}=0$, and so
$v=0$. Hence the state trajectory $z(\cdot)$ of $\Sigma^\ka$ is given
by $z(t)=\tline_t^\ka x$. Since the systems $\Sigma$ and $\Sigma^\ka$
have the same state trajectory $z(\cdot)$, we obtain \m
$z(t)=\tline_t^\ka x=\tline_t\m x$ for all $t\geq 0$. In particular,
it follows that $\|\tline_t x\|=\|x\|$ for all $t\geq 0$. We now redo
this whole argument with the dual of \m $\Sigma$ using Proposition
\ref{X}. This shows that $\|\tline_t^*x\|=\|x\|=\|\tline_tx\|$ for all
$t\geq 0$.

The output signal $y$ of $\Sigma$ with initial state $x\in X^u$ and
with zero input equals the output signal of $\Sigma^\ka$. We have seen
earlier that $\hat{y}=0$, i.e., $x\in\Nscr$. According to \rfb{Cweak}
we obtain $x=0$. Hence $X^u=\{0\}$, so that $\tline^\ka$ is weakly
stable. Dual arguments show that \rfb{Bweak} implies the weak
stability of $\tline^{\ka*}$, which in turn is equivalent to the weak
stability of $\tline^\ka$.

(b) The second step is to derive a useful estimate concerning the
open-loop transfer function $\GGG$. Since $cI+\GGG$ is a positive 
transfer function (see Remark \ref{Aug_2_2007} and Corollary
\ref{cpassive}), for all $s\in\cline_0$ and for all $u_0\in U$ with
$\|u_0\|=1$ there holds
$$ \left\| \left( \frac{1}{\ka} I + \GGG(s) \right) u_0 \right\| 
   \m\geq\m \Re \left\langle \left( \frac{1}{\ka} I + \GGG(s) \right)
   u_0,u_0 \right\rangle \qquad\m$$
\BEQ{eig2}
   \m\qquad =\m \frac{1}{\ka} -c \mm+\mm
   \Re \left\langle \left( cI + \GGG(s) \right) u_0,u_0 \right\rangle
   \m\geq\m \frac{1}{\ka} -c \m>\m 0 \m.
\EEQ
Suppose that $\o\in\rline$ is such that $\jw_0\in\rho(A)$. Then
$(\frac{1}{\ka}I+\GGG(s))u_0$ is analytic on a small neighbourhood
$N_\e=\{z|\ \|z-\jw_0\|<\e\}$ (and of course also on $\cline_0$). By
continuous extension we obtain from \rfb{eig2} that for $\jw_0\in
\rho(A)$ and $\|u_0\|=1$,
\BEQ{eig3}
   \left\| \left( \frac{1}{\ka} I + \GGG(\jw_0) \right) u_0 \right\| 
   \m\geq\m \frac{1}{\ka} -c \m>\m 0 \m.
\EEQ

(c) To prove strong stability, first we show that 
\BEQ{countable}
   \sigma(A^\ka)\cap i\rline \subset  \sigma(A)\cap i\rline \m.
\EEQ
Since $-\kappa I$ is an admissible feedback operator for $\Sigma$ that
produces the well-posed closed-loop system $\Sigma^\kappa$, we can
apply Theorem \ref{pertnode} to obtain that if $\l\in\sigma(A^\ka)\cap
\rho(A)$, then $-\frac{1}{\ka}\in\sigma(\GGG(\l))$. In particular,
suppose that for some $\o\in\rline$ we have $\jw_0\in\sigma(A^\ka)\cap
\rho(A)$, then $-\frac{1}{\ka}\in\sigma(\GGG(\jw_0))$. To prove \rfb
{countable}, we have to show that such an $\o_0$ cannot exist. We 
consider the classical three cases:

If $-\frac{1}{\ka}$ is in the point spectrum of $\GGG(\jw_0)$, then
there exists $u_0\in U$ with $\|u_0\|=1$ such that
$$\frac{1}{\ka} u_0 + \GGG(i\o_0)u_0 \m=\m 0 \m.$$
This is in direct contradiction with \rfb{eig3}.

If $-\frac{1}{\ka}$ is in the residual spectrum of $\GGG(\jw_0)$, then
it is in the point spectrum of $\GGG(i\o_0)^*$. Then a similar
argument applied to the dual system node with the dual transfer
function $\GGG^d(s)=\GGG(\bar{s})^*$ also leads to a contradiction.

The remaining possibility is that $-\frac{1}{\ka}$ is in the
continuous spectrum of $\GGG(\jw_0)$. In this case there exists a
sequence $u_n\in U$ with $\|u_n\|=1$ and
\BEQ{eign}
    \left\| \frac{1}{\ka} u_n + \GGG(i\o_0)u_n \right\|
    \rarrow\, 0 \mbox{~~~as~~} n\m\rarrow\, \infty \m.
\EEQ
This is again in contradiction with \rfb{eig3}. We have established
\rfb{countable}.

This implies that (with the assumption in the last part of the
theorem) $\sigma(A^\ka)\cap i\rline$ is countable. According to
Theorem \ref{diss} $\tline^\ka$ is a contraction semigroup. Its weak
stability (shown earlier in this proof) implies that $A^\ka$ and
$A^{\ka*}$ have no eigenvalues on $i\rline$. So we can apply the main
result from Arendt and Batty \cite{AreBat} to establish the strong
stability of $\tline^\ka$. Since $\sigma({A^\ka}^*)\cap i\rline=\mm
\overline{\sigma(A^\ka)\cap i\rline}$, a similar argument shows that
$\tline^*$ is also strongly stable. \kasten

\begin{remark} \label{weakobs} {\rm
Note that \rfb{Cweak} means that the unitary part of $\tline$ (i.e.,
its restriction to the space $X^u$ defined in the last proof) is
approximately observable in infinite time. In particular, \rfb{Cweak}
holds if $(A,C)$ is approximately observable in infinite time. Similar
remarks hold for the dual condition \rfb{Bweak}. In particular,
\rfb{Bweak} is implied by the approximate controllability in infinite
time of $(A,B)$.}
\end{remark}

\begin{remark} {\rm
Theorem \ref{stability} generalizes Theorem 14 in Batty and Phong
\cite{BatPho} to unbounded operators $B$ and $C$, while eliminating
the assumption that $C=B^*$.}
\end{remark}

\begin{remark} {\rm
Let $\Sigma_i$, $i=1,2$, be two impedance passive system nodes with
the inputs $u_i$, outputs $y_i$ and transfer functions $\GGG_i$,
$i=1,2$. Hence we have
$$  \|z_i(\tau)\|^2-\|z_i(0)\|^2 \leq 2 \int_0^\tau {\rm Re}\langle 
    u_i(t),  y_i(t) \rangle \dd t \m.$$
Subsituting $u_1=-y_2+v_{cl},\ \ u_2=y_1$ we obtain
$$  \|z_1(\tau)\|^2-\|z_1(0)\|^2\leq 2\int_0^\tau {\rm Re}\langle 
    -y_2(t)+v_{c}(t),y_1(t) \rangle \dd t \m,$$
$$  \|z_2(\tau)\|^2-\|z_2(0)\|^2\leq 2\int_0^\tau {\rm Re}\langle 
    y_1(t),y_2(t) \rangle \dd t \m,$$
and adding these two equations yields
$$  \|z_1(\tau)\|^2+\|z_2(\tau)\|^2-\|z_1(0)\|^2-\|z_2(0)\|^2
    \leq 2\int_0^\tau {\rm Re}\langle v_{c}(t),y_1(t) \rangle 
    \dd t \m.$$
Hence the above feedback connection produces another impedance passive
system $\Sigma_{cl}$ with the state $\bbm{ z_1\\z_2}$, the input
$v_{cl}$, the output $y_{cl}=y_1$ and the transfer function $\GGG_{cl}
\mm=\GGG_1(I-\GGG_2\GGG_1)^{-1}$ (see van der Schaft \cite[Proposition
3.2.5]{Schaft}).

Hence, according to Theorem \ref{diss}, for all $\e>0$ the feedback
control law $v_{cl}=-\e y_{cl}+v$ produces a well-posed closed-loop
system $\Sigma_\e$ that is system stable. Its transfer function is
given by \m $\GGG_\e=\GGG_1(I+\GGG_2\GGG_1+\e\GGG_1)^{-1}$. Note that
$\Sigma_\e$ can equally well be obtained by applying the
transformation $u_1=-(\e I+\GGG_2)y_1$ to $\Sigma_1$. In the
finite-dimensional literature this is known as dynamic feedback
stabilization.}
\end{remark}

\section{A class of damped second order systems with colocated
         actuators and sensors} 

\ \ \ In this section we introduce a class of damped second order
systems that are impedance passive system nodes (in the sense that
\rfb{statepassive} holds), where the actuators and sensors are
colocated (in the sense that assumption {\bf COL} holds).

Let $U_0$ and $H$ be Hilbert spaces and let $A_0:\Dscr(A_0) \rarrow\m
H$ be positive and boundedly invertible on $H$. For every $\mu>0$, we
define $H_\mu=\Dscr(A_0^\mu)$, with the norm $\|\varphi\|_\mu=\|A_0
^\mu\varphi\|_H$, and we define $H_{-\mu}=H^*_\mu$ (duality with
respect to the pivot space $H$). We denote $H_0=H$ and $\|\varphi
\|_0=\|\varphi\|_H$. We assume that
$$ C_0 \m\in\m \Lscr(H_\half,U_0) \m,\quad M \m\in\m \Lscr(H_\half,
   H_{-\half}) \m,\ \ \ \ M\geq 0 \m.$$
By $M\geq 0$ we mean that it defines a positive quadratic form on
$H_\half$. We identify $U_0$ with its dual, so that
$C_0^*\in\Lscr (U_0,H_{-\half})$.

We wish to study the system described by the abstract second order
equation 
\BEQ{colfirst} 
   \ddot{q} + M \dot q + A_0 q \m=\m C_0^* u \m,\qquad
   y \m=\m C_0 \m \dot{q}.
\EEQ
This is a slight generalization of the class of systems studied in
Tucsnak and Weiss \cite{TuWe03,WeTu}. We introduce the state space
$X=H_\half\times H$ and we define
\BEQ{A}
   A:\Dscr(A) \rarrow X \m,\ \ \ A \m=\m \bbm{ 0 & I \\ -A_0 & -M},
\EEQ
$$ X_1 \m=\m \Dscr(A) \m=\m \left\{ \sbm{q\\ w} \in H_\half \times
   H_\half \ |\ A_0 q + M w \in H \right\} \m,$$
$$ C \m=\m \left[\begin{array}{cc} 0 & C_0\end{array}\right] \m,\ \ \ 
   B \m=\m \bbm{ 0 \\ C_0^* } \m,$$
then $B\in\Lscr(U_0,X_{-1})$, $B^*,C\in\Lscr(X_1,U_0)$. It is also
easy to verify that {\bf COL} holds and 
$$ A^*:\Dscr(A^*) \rarrow X \m,\ \ \ A^* \m=\m \bbm{ 0 & -I \\
   A_0 & -M},$$
$$ \Dscr(A^*) \m=\m \left\{ \sbm{q\\ w} \in H_\half \times
   H_\half \ |\ -A_0 q + M w \in H \right\} \m,$$
It can be shown that the operators $A,B,C$ and $D$ define a compatible
system node $\Sigma$ on $(U,_0X,U_0)$, whose state trajectories $z$
and output functions $y$ corresponding to $C^2$ input functions $u$
and compatible initial conditions satisfy
\begin{eqnarray} \label{colsecond}
   \left\{ \begin{array}{rcl} \dot{z} &=& A z + B u \m,\\
   y &=& \overline{C} z \m. \end{array} \right.
\end{eqnarray}

As in \cite{WeTu}, \cite{CuWeABB*}, it is not difficult to prove that
\rfb{colfirst} is equivalent to \rfb{colsecond}, for any $C^2$ input
signal $u$ and compatible initial conditions $q(0)$ and $\dot{q}(0)$,
i.e., $A_0q(0)+M\dot q(0)-C_0^*u(0)\in H$. Indeed, the connection is
provided by
$$z \m=\m \bbm{q\\ \dot q} \m.$$
Note that
$$A+ A^* \m=\m \bbm{ 0 & -0 \\0 & -2M},$$
and so {\bf ESAD} is not satisfied if $M$ is unbounded. However, since
$0\in \rho(A)$, it is relatively easy to show that this system is
impedance passive by using Corollary \ref{corASS}. We have
$$ A^{-1} \m= \bbm{-A_0^{-1} M & -A_0^{-1}\\ I&0},\ A^{-*} \m=
   \bbm{-A_0^{-1} M & A_0^{-1}\\ -I&0},$$
and so
$$ CA^{-1} \m= \bbm{ C_0 & 0}\m =\m- B^*A^{-*} \m.$$ 
Moreover, $\GGG(0)=0$ and the conditions of Corollary \ref{corASS} 
are satisfied. 

{\em Example.} In the literature there are many examples of systems of
the type discussed above (see Luo {\em et al} \cite{LuoGuoMor98}). We
consider an example in Bontsema \cite{Bontsema}. This is an idealized
model of a large flexible satellite with a central hub. The partial
differential equation model for a beam of length 2 is
$$ \rho a \frac{\partial^2w}{\partial t^2}(t,x)+\bar{E}I 
   \frac{\partial^5w}{\partial t\partial x^4}(t,x)+
   EI \frac{\partial^4w}{\partial x^4}(t,x) \m=\m 0,$$
$$ \frac{\partial^3w}{\partial x^3}(-1,t)\m=\m 0 \m=\m
   \frac{\partial^3w} {\partial x^3}(1,t) \m,\ \ \frac{\partial^2w}
   {\partial x^2}(-1,t)\m=\m 0 \m=\m
   \frac{\partial^2w}{\partial x^2}(1,t),$$
$$ u_0(t) \m=\m EI[\frac{\partial^3w}{\partial x^3}(0^+,t) -\frac
   {\partial^3w}{\partial x^3}(0^-,t)]
 +\bar{E}I[\frac{\partial^4w}{\partial t\partial x^3}(0^+,t)-\frac
   {\partial^4w}{\partial t\partial x^3}(0^-,t)],$$
$$ u_1(t) \m=\m -EI[\frac{\partial^2w}{\partial x^2}(0^+,t)- \frac
   {\partial^2w}{\partial x^2}(0^-,t)]+\bar{E}I[\frac{\partial^3w}
   {\partial t\partial x^2}(0^+,t)-
   \frac{\partial^3w}{\partial t\partial x^2}(0^-,t)],$$
$$ y_0(t) \m=\m \frac{\partial w}{\partial t}(0,t),\ \ y_1(t) \m=\m
   \frac{\partial^2 w}{\partial t\partial x}(0,t).$$
Here, $w(x,t)$ represents the vertical displacement of the beam at
co-ordinate $x$ along the beam at time $t$, $a$ is the cross-sectional
area of the beam, $\rho$ its mass density, $E$ is Young's modulus, $I$ the moment of inertia of the beam per cross-section and $\bar{E}$ is a constant reflecting the stress-strain relation in the beam. $u_0(t)$
is the force and $u_1(t)$ is the moment acting on the centre of the
beam at time $t$. $y_0(t)$ and $y_1(t)$ are the measurements of the
velocity, respectively the angular velocity, in the middle of the beam
at time $t$. Although the observation and control operators are
admissible, the transfer function grows as $\sqrt{s}$ as $S\to \infty$ along the real axis and so it is not well-posed. Various stabilization schemes for this unstable system were developed in \cite{Bontsema}, but stabilization using colocated actuators and sensors was not treated. Here show that it can be strongly stabilized using static output feedback.
Denote by $H^4(-1,1)$ the Sobolev space
$$H^4(-1,1) \m=\m \{f\in L_2(-1,1)\ | \ \frac{df}{dx}, \frac{df^2}{dx^2}\frac{d^3f}{dx^3} \in L_2(-1,1)\}, $$
where the derivatives are defined in terms of distributions.
In \cite{Bontsema} it was shown that 
the operator $A_0=\frac {d^4}{dx^4}$ with domain 
$$\Dscr(A_0) \m=\m \{ w\in H^4(-1,1) \ | \ \frac{d^2w}{dx^2}(-1)=0, \ \frac{d^3w}{dx^3}(-1)=0, \ \frac{d^2w}{dx^2}(1)=0, \ \frac{d^3w}{dx^3}(1)=0\}$$ 
is a densely defined, nonnegative, self-adjoint operator on $L_2(-1,1)$. 
Since it has a double eigenvalue at zero, to fit this example into 
the framework of \rfb{A} with $H:=L_2(-1,1)$, $M=A_0$,
we need to introduce the modified inner product on the state-space $X=\Dscr(A_0)\oplus L_2(-1,1)$ 
$$ \langle \sbm{z_1\\z_2},\sbm{w_1\\w_2}\rangle_X \m=\m \langle z_1,
   w_1 \rangle +\langle A_0^{\half} z_1,A_0^{\half}w_1\rangle +
   \langle z_2,w_2\rangle,$$
where the inner products are in $L_2(-1,1)$. 
Let $U_0=\cline^2$,  $y=\bbm{y_0&y_1}$ , $u=\bbm{u_0\\u_1}$ and define
$$C_0 w \m=\m \bbm{\frac{\partial w}{\partial t}(0)&\frac{\partial^2 w}{\partial t\partial x}(0)}.$$
Although this example now fits into the framework 
\rfb{A}--\rfb{colfirst}, due to the eigenvalues at zero, Corollary
\ref{ASS} is not applicable. Instead we use Theorem \ref{tpassive} to
show that it is impedance passive with the colocated pair $u,y$.  First
we compute
$$(sI-A)^{-1}=\m=\m \bbm{ V(s)(sI+M) & V(s) \\-V(s)A_0 & sV(s)},$$
where $V(s)=(s^2I+Ms+A_0)^{-1}\in \Lscr(H_{-\half},H_{\half})$\m.
$$ (sI-A)^{-1}B \m=\m \bbm{ V(s)C_0^*  \\sV(s)C_0^*},$$
$$(sI+A^*)(sI-A)^{-1}B-B^*\m=\m(A+A^*)(sI-A)^{-1}B\m,$$
$$\GGG(s)\m=\m C(sI-A)^{-1}B\m=\m B^*sV(s)B \m.$$
 From here, \vspace{-1mm}
$$ \GGG(s)+\GGG(s)^*-B^*(\overline{s}I-A)^{-1}2{\rm Re}s
   (\overline{s}I-A)^{-1}B$$
\vspace{-7mm}
\begin{eqnarray*}
   &=&\m -B^*V(\overline{s})2|s|^2MV(s)B\\
   &=&\m B^*(\overline{s}I-A^*)^{-1}(A+A^*)(sI-A)^{-1}B\m.
\end{eqnarray*}
As in the proof of Proposition \ref{job} the inequality in part (a) of
Theorem \ref{tpassive} can be factored as \rfb{matrixineq} and the
system $\Sigma$ is impedance passive. This example is unstable due to
a double eigenvalue at $0$. In \cite{Bontsema} it is shown that the
spectrum of its generator is countable, it is approximately
controllable in infinite time and approximately observable in infinite
time. We have seen that it is impedance passive, hence for any
$\ka>0$ the feedback $u=-\ka y+v$ results in a strongly stable
closed-loop system.

\section{The effect of damping and feedthrough on a class of
         second order systems} 

\ \ \ In this section we examine the conditions under which the
following second order system will be almost impedance passive.
$$ \ddot{q} + M \dot q + A_0 q \m=\m B_0 u  \m,$$
$$ y \m=\m C_0 \m \dot{q},$$
where $A_0,M,C_0$ are as in Section 6 and 
$B_0\in\Lscr(U_0,H_{-\half})$. As in Section 6 it can be shown that
this is equivalent to a compatible system node with everything 
defined as in Section 6, except for $B$ which is given by
$$ B \m=\m \bbm{ 0 \\ B_0 } \m, \mbox{~~and}\ \ B^* \m=\m
   \bbm{ 0 & B_0^* } \m.$$
As before, $\GGG(0)=0$. We apply Proposition \ref{ABCimp} with
$\omega=0$ to obtain conditions for $\Sigma_E$ to be impedance
passive. We compute 
$$ A^{-1}+A^{-*}\m=\bbm{-2A_0^{-1}M&0\\0&0},$$
$$ CA^{-1}+B^*A^{-*}=\bbm{C_0-B_0^*&0},\mbox{~~and}$$
$$ A^{-*}C^*+A^{-1}B=\bbm{A_0^{-1}(C_0^*-B_0)\\0}.$$
So $\Sigma_E$ is impedance passive if and only if
$$ 2 {\rm Re}\langle u,(C_0-B_0^*)x\rangle  \leq 2\langle Mx,x
   \rangle-\langle u,Eu\rangle$$
for all $x\in H_{\half}$, $u\in U_0$. For simplicity suppose that
$M$ is invertible. Then by completing the square you can show that
with $E=-\frac{1}{4} (C_0-B_0^*)M^{-1}(C_0^*-B_0)$, $\Sigma_E$ is
impedance passive. So if $M$ is invertible, $\Sigma$ is almost
impedance passive for an arbitrary choice of $C_0-B_0^*$. The
actuators and sensors need not be colocated. Note that colocated
actuators and sensors are necessary for $\Sigma$ to be impedance
passive.

\section{A class of damped second order systems with colocated
        actuators and sensors, but not impedance passive} 

\ \ \ In this section we show that colocated actuators and sensors
need not imply impedance passivity. We then design a non-colocated
input and output pair so that the system is almost impedance passive
and we derive an explicit expression for the minimal $E$ for which
\rfb{relaxed} holds.  Following Section 6, we formulate our systems
as compatible system nodes. Let $U_0, U_1, H$, be Hilbert spaces and
$H_\mu=\Dscr(A_0^\mu)$, $H_{-\mu}=H^*_\mu$ and $M$ be as in Section 6.
We assume that
$$ C_0 \m\in\m \Lscr(H_\half,U_0) \m,\qquad
   C_1 \m\in\m \Lscr(H_1,U_1) \m.$$
As before, we identify $U_0$ and $U_1$ with their duals, so that
$C_0^*\in\Lscr (U_0,H_{-\half})$, \m $C_1^*\in\Lscr(U_1,H_{-1})$.

We assume that $C_0$ and $C_1$ have extensions $\overline{C_0}$ and
$\overline{C_1}$ such that the operators
$$ D_0 \m=\m \overline{C_0}A_0^{-1}C_1^*\m\in\m\Lscr(U_1,U_0)\m,$$
$$ D_1 \m=\m \overline{C_1}A_0^{-1}C_0^*\m\in\m\Lscr(U_0,U_1)$$
$$ C_2 \m=\m \overline{C_1}A_0^{-1}M\m\in\m\Lscr(H_\half,U_1)$$
exist. Moreover, we assume that $C_2$ has an extension $\overline
{C_2}$ such that the following operator exists
$$ D_2 \m=\m \overline{C_2}A_0^{-1}C_1^*\m\in\m\Lscr(U_1,U_1).$$
We define again $X=H_\half\times H$ and $A$ is defined by \rfb{A},
$$ C \m=\m \left[\begin{array}{cc} 0 & C_0, \\
   \overline{C_1} & 2C_2 \end{array}\right] \m,\ \ \ 
   \overline{C} \m=\m \left[\begin{array}{cc} 0 & \overline{C_0},\\
   \overline{C_1} & 2\overline{C_2} \end{array}\right] \m,$$
$$ B \m=\m \bbm{ 0 & A_0^{-1} C_1^* \\ C_0^* & 0},\ \ \ 
   D \m=\m \left[\begin{array}{cc} 0 & D_0 \\
   0 & D_2 \end{array}\right] \m,$$
   $$D_{col}\m=\m \bbm{0&D_0\\ 0&0} \m,$$
where $B\in\Lscr(U,X_{-1})$, $B^*,C\in\Lscr(X_1,U)$. It is 
easy to verify that
$$ B^* \m=\m \bbm{ 0 & C_0 \\ \overline{C_1} & 0}$$
and it can be shown that the operators $A,B,\overline{C}$ and $D$
and the operators $A,B,B^*$ and $D_{col}$ define compatible system
nodes $\Sigma$ and $\Sigma_{col}$, respectively, on $(U,X,U)$. We
denote their transfer functions by $\GGG$, respectively,
$\GGG_{col}$. The state trajectories $z$ and output functions $y$ of
$\Sigma$ corresponding to $C^2$ input functions $u$ and compatible
initial conditions satisfy
\begin{eqnarray} \label{second}
   \left\{ \begin{array}{rcl} \dot{z} &=& A z + B u \m,\\
   y &=& \overline{C} z + D u \m. \end{array} \right.
\end{eqnarray}
The state trajectories $z$ and output functions $y$ of $\Sigma_{col}$
corresponding to $C^2$ input functions $u$ and compatible initial
conditions satisfy
\begin{eqnarray} \label{third}
   \left\{ \begin{array}{rcl} \dot{z} &=& A z + B u \m,\\
   y &=& \overline{B^*} z + D_{col} u \m. \end{array} \right.
\end{eqnarray}

As in \cite{CuWeABB*} it is straightforward to prove that \rfb{second}
is equivalent to the following second order differential equation and
the two output equations:
\BEQ{first} 
   \ddot{q} + M \dot q + A_0 q \m=\m C_0^* u_0 +
   A_0^{-1} C_1^* \dot u_1 \m,
\EEQ
\BEQ{first-prime}
   y_0 \m=\m \overline{C_0} \m \dot{q},\ \ \ \ y_1 \,=\, 
   \overline{C_1} \m q +\overline{C_2} \m\dot{q}
\EEQ
for $C^2$ input signals $u_0,u_1$ and compatible initial
conditions $q(0)$ and $\dot{q}(0)$. Indeed, the connection is
provided by 
$$ z \m= \bbm{q\\ w}\m,\ \ u \m= \bbm{u_0\\ u_1}\m,
   \ \ y \m= \bbm{y_0\\ y_1} \m,\ \ w=\dot q-A_0^{-1}C_1^* u_1 \m.$$
\medskip
Similarly, \rfb{third} is equivalent to the same second order 
differential equation \rfb{first} and two output equations:
\BEQ{third-prime}
   y_0 \m=\m \overline{C_0} \m \dot{q},\ \ \ \ y_1 \,=\, 
   \overline{C_1} \m q  \m\
\EEQ
for $C^2$ input signals $u_0,u_1$ and compatible initial
conditions $q(0)$ and $\dot{q}(0)$. Applying Proposition
\ref{ABCimp} we obtain the following necessary and sufficient
conditions for $\Sigma_{col}$ to be impedance passive:
$$ \langle \overline{u},(\GGG_{col}(0)+\GGG(0)_{col}^*)\overline{u}
   \rangle\geq 2\langle x,Mx\rangle +2{\rm Re}\langle C_2x,u_1
   \rangle$$
for all $x\in H_{\half}, u_0\in U_0, u_1\in U_1$, where
$\overline{u}=\bbm{u_0\\u_1}$. But
$$\GGG_{col}(0)+\GGG_{col}(0)^*=\bbm{0&D_1^*\\D_1&D_2+D_2^*}$$
and so $\Sigma_{col}$ is not impedance passive. It will be almost
impedance passive with 
$$2E=\bbm{0&-D_1^*\\-D_1&2E_{22}-D_2-D_2^*},$$
provided that
$$2\langle u_1,E_{22}u_1\rangle\geq
2\langle x,Mx\rangle +2{\rm Re}\langle C_2x,u_1\rangle.$$
for all $x\in H_{\half}, u_1\in U_1$. If $M$ is invertible, then this
is satisfied with $E_{22}\geq-\frac{1}{4}C_2M^{-1}C_2^*$.
Alternatively, if there exists a $Q\in\Lscr(X,U_1)$ such that
$QM^{\half}z=\half C_2z$ for $x\in H_{\half}$, then $\Sigma _E$ is
impedance passive with $E_{22}\geq -Q^*Q$.  Otherwise, it is hard to
see how to choose $E_{11}$.  Instead we show that by choosing the
output as in \rfb{first-prime}, it is clear how to choose an $E$ so
that $\Sigma_E$ is impedance passive. By applying Corollary
\ref{corASS}, we obtain sufficient conditions for $\Sigma_E$ to be
impedance passive.

\begin{lemma}\label{example}
Let $\Sigma$ be the compatible system node from \rfb{second} and
denote its transfer function by $\GGG$. Then for $E=E^*\in \Lscr(U)$,
$\Sigma_E$ will be impedance passive if and only if
$$\GGG(0)+\GGG(0)^*+2E \m\geq\m 0 \m.$$
Moreover, the smallest $E$ for which this holds is 
\BEQ{CEeqn}
   E \m=\m -\half \left[ \begin{array}{cc} 0&D_1^* \\
   D_1&0 \end{array}\right] \m.
\EEQ
\end{lemma}

\medskip
{\it Proof .} \m First we show that condition \rfb{ASS} holds.
We have $0\in\rho(A)$ and 
$$ CA^{-1} \m= \bbm{ C_0 & 0\\ C_2&
   -C_1 A_0^{-1} }$$
and
$$ B^*A^{-*} \m= \bbm{ -C_0 & 0 \\ -C_2 & C_1 A_0^{-1}}.$$ 
So we have $CA^{-1}+B^*A^{-*}=0$, and applying Corollary \ref{corASS}
to $\Sigma_E$ proves that it is impedance passive. 

The transfer function of $\Sigma$ is given by
$\GGG(s)=\overline{C}(sI-A)^{-1}B+D$, which is easy to compute in
terms of $A_0,C_0$ and $C_1$. We have
\BEQ{G_beam}
   \overline{C}A^{-1}B \m=\m 
   \left[ \begin{array}{cc} 0&\overline{C_0}A_0^{-1}C_1^*  \\
   -\overline{C_1}A_0^{-1}C_0^*&\overline{C_2}A_0^{-1}C_1^*
   \end{array} \right] \m,
\EEQ
and 
$$ \GGG(0)+\GGG(0)^* \m= \left[ \begin{array}{cc} 0&D_1^* \\ 
   D_1&0 \end{array} \right] \m.$$
So the smallest $E$ such that $\Sigma_E$ is impedance passive is
given by \rfb{CEeqn}. \kasten

An example of this type of system with no damping was given in Weiss
and Curtain \cite{WCbeam}. It was a model of a hinged elastic beam with
2 sensors, one measuring the curvature and one measuring the angular
velocity at a point on the beam.
 


\begin{thebibliography}{}
\small

\bibitem{AmLiTu} K. Ammari, Z. Liu and M. Tucsnak, \m Decay rates
 for a beam with pointwise force and moment feedback, \m {\em
 Mathematics of Control, Signals and Systems} {\bf 15} (2002),
 pp.\,229--255.

\bibitem{AreBat} W. Arendt and C.J.K. Batty, \m Tauberian theorems and
 stability of one-parameter semigroups, \m {\it Transactions of the
 American Mathematical Society} {\bf 306} (1988), pp.\,837--841.

\bibitem{ABHN} W. Arendt, C.J.K. Batty, M. Hieber and F. Neubrander.
 \m {\em Vector-valued Laplace Transforms and Cauchy Problems}, \m
 Birkh\"{a}user Verlag, Basel, 2001.

\bibitem{BaiHub} T. Bailey and J.E. Hubbard jr., \m Distributed
 piezoelectric polymer active vibration control of a cantilever beam,
 \m {\it AIAA Journal on Guidance, Control and Dynamics} {\bf 8}
 (1985), pp.\,605--611.

\bibitem{Bal1} A.V. Balakrishnan, \m Compensator design for stability
 enhancement with collocated controllers, \m {\it IEEE Trans. Autom.
 Control} {\bf 36} (1992), pp.\,994--1008.

\bibitem{Bal2} A.V. Balakrishnan, \m Shape control of plates with
 piezo actuators and collocated position\m/\m rate sensors, \m {\it
 Applied Math. and Comput.} {\bf 63} (1994), pp.\,213--234.

\bibitem{BatPho} C.J.K. Batty and V.Q. Phong, \m Stability of
 individual elements under one-parameter semigroups, \m {\it Trans.
 Amer. Math. Soc.} {\bf 322} (1990), pp.\,805--818.

\bibitem{Ben} C.D. Benchimol, \m A note on weak stabilizability of
 contraction semigroups. \m {\it SIAM Journal on Control and Optim.}
 {\bf 16} (1978), pp.\,373--379.

\bibitem{Bontsema} J. Bontsema, \m Dynamic Stabilization of Large
 Flexible Space Structures, Ph.D.Thesis, Rijksuniversiteit Groningen,
 The Netherlands, 1989.

\bibitem{CurWei} R.F. Curtain and G. Weiss, \m Well-posedness of 
 triples of operators (in the sense of linear systems theory), \m 
 {\it Control and Estimation of Distributed Parameter Systems} 
 (F. Kappel, K. Kunisch, W. Schappacher, eds.), pp.\,41--59,
 Birkh\"auser-Verlag, Basel, 1989. 

\bibitem{CuWeABB*} R.F. Curtain and G. Weiss, \m Exponential
 stabilization of well-posed systems by colocated feedback, \m 
 {\it SIAM J. Control and Optim.}, {\bf 45} (2006), pp.\,273--297.

\bibitem{CurZwa95} R.F. Curtain and H.J. Zwart, \m {\it An 
 Introduction to Infinite-Dimensional Linear Systems Theory}, 
 Springer-Verlag, New York, 1995. 

\bibitem{CurJac} R.F. Curtain and B. Jacob, \m Spectral properties
 of pseudo-resolvents under structured perturbations, \m {\it
 Mathematics of Control, Signals and Systems}, {\bf 21} (2008), 
 pp.\,21--50.

\bibitem{Davi} E.B. Davies, \m {\it One-Parameter Semigroups}, \m
 Academic Press, London, 1980.

\bibitem{EnNa} K. Engel and R. Nagel, \m {\it One-parameter
 Semigroups for Linear Evolution Equations}, \m Graduate Texts in 
 Math. vol. 194, Springer-Verlag, New York, 2000.


\bibitem{Gib} J.S. Gibson, \m A note on stabilization of infinite 
 dimensional linear oscillators by compact feedback, \m {\it SIAM J. 
 Control and Optim.} {\bf 18} (1980), pp.\,311--316.

\bibitem{GorZwaMas05} Y.~Le Gorrec, H.J.~Zwart and B.~Maschke, \m 
 Dirac structures and boundary control systems associated with 
 skew-symmetric differential operators. {\em SIAM J. of Control and 
 Optim.}, {\bf 44} (2005), pp.\,1864--1892.

\bibitem{GuoLuo} B-Z. Guo and Z-H. Luo, \m Controllability and 
 stability of a second-order hyperbolic system with colocated 
 sensor/actuator, \m {\it Systems \& Control Letters} {\bf 46}
 (2002), pp.\,45--65.


\bibitem{Haraux} A. Haraux, \m Une remarque sur la stabilisation de
 certains syst\`emes du deuxi\`eme ordre en temps, \m {\it 
 Portugaliae Mathematica} {\bf 46} (1989), pp.\,245--258.

\bibitem{LasTrig} I. Lasiecka and R. Triggiani, $L_2(\Sigma)$-regular%
 ity of the boundary to boundary operator $B^*L$ for hyperbolic and
 Petrowski PDE's, \m {\it Abstract and Applied Analysis} {\bf 19} 
 (2003), pp.\,1061--1139. 
   
\bibitem{LasTrig3} I. Lasiecka and R. Triggiani, \m {\it Control
 Theory for Partial Differential Equations: Continuous and
 Approximation Theories. {II}, Abstract hyperbolic-type systems
 over a finite time horizon}, Encyclopedia of Mathematics and its
 Applications {\bf 75}, Cambridge University Press, Cambridge, 2000.
   
\bibitem{Liu97} K. Liu, \m Local distributed control and damping for 
 the conservative systems, \m {\it SIAM J. Control and Optim.} {\bf 35}
 (1997), pp.\,1574--1590. 

\bibitem{LivsBook} M.S. {Liv\v sic}, {\it Operators, Oscillations,
 Waves (Open Systems)}, volume~34 of {\em Translations of Mathematical
 Monographs}. American Mathematical Society, Providence, Rhode Island,
 1973.

\bibitem{LyPh} Y.I. Lyubich and V.Q. Phong, \m Asymptotic stability of
 linear differential equations in Banach spaces, \m {\em Studia Math.}
 {\bf 88} (1988), pp.\,37--42.  

\bibitem{LuoGuoMor98} Z-H. Luo, B-Z. Guo and O.Morgul, \m {\it 
 Stability and Stabilization of Infinite Dimensional Systems with 
 Applications}, \m Springer-Verlag, London, 1999.

\bibitem{MalStafWei} J.Malinen, O.J. Staffans and G. Weiss, \m When
 is a linear system conservative? \m {\it Quarterly of Applied
 Math.}, to appear in 2006.

\bibitem{NagFoi} B. Sz.-Nagy and C. Foias, \m {\it Harmonic Analysis
 of Operators on Hilbert Space}, \m North-Holland, Amsterdam, 1970
 (transl. of the French edition of 1967).

\bibitem{Oos} J.C. Oostveen, \m {\em Strongly Stabilizable 
 Infinite-Dimensional Systems}, \m Frontiers in Applied Mathematics,
 SIAM, Philadelphia, 2000.
 
\bibitem{Opm} M.R. Opmeer, \m Infinite-dimensional linear systems:
 a distributional approach, {\it Proc. London Math. Society} {\bf 91}
 (2005), pp.\,738--760.
 
\bibitem{RoRo} M. Rosenblum and J. Rovnyak, \m {\em Hardy Classes
 and Operator Theory}, \m Oxford University Press, 1985.

\bibitem{Rudin} W. Rudin, \m {\em Real and Complex Analysis}, \m
 McGraw-Hill, New York, 1966.

\bibitem{Rus1} D.L. Russell, \m Linear stabilization of the linear
 oscillator in Hilbert space, \m {\it J. Math. and Applications}
 {\bf 25} (1969), pp.\,663--675.

\bibitem{sala} D. Salamon, \m Infinite dimensional systems with
 unbounded control and observation: A functional analytic approach,
 \m {\it Trans. Amer. Math. Soc.} {\bf 300} (1987), pp.\,383--431.

 
\bibitem{Schaft} A. van der Schaft, \m {\em $L_2$-Gain and Passivity
 Techniques in Nonlinear Control}, second (enlarged) edition, \m
 Springer-Verlag, London, 2000.

\bibitem{Slem1} M. Slemrod, \m A note on complete controllability
 and stabilizability for linear control systems in Hilbert space, \m
 {\it SIAM J. Control and Optim.} {\bf 12} (1974), pp.\,500--508.
 
\bibitem{Slem2} M. Slemrod, \m Stabilization of boundary control
 systems, \m {\it J. of Diff. Equations} {\bf 22} (1976),
 pp.\,402--415.

\bibitem{Slem3} M. Slemrod, \m Feedback stabilization of a linear
 control system in {Hilbert} space with an a priori bounded control,
 \m {\it Math. of Control, Signals and Systems} {\bf 2} (1989),
 pp.\,265--285.

\bibitem{stafTAMS} O.J. Staffans, \m Quadratic optimal control of
 stable well-posed linear systems, \m {\it Trans. Amer. Math.
 Society} {\bf 349} (1997), pp.\,3679--3715.


\bibitem{StafPas} O.J. Staffans, \m Passive and conservative 
 continuous-time impedance and scattering systems. Part I: Well-posed
 systems, \m {\it Mathematics of Control, Signals and Systems}
 {\bf 15} (2002), pp.\,291--315.

\bibitem{StafCol} O.J. Staffans, \m Stabilization by collocated
 feedback, \m {\it Directions in Mathematical Systems Theory and
 Optimization}, A. Rantzer and C.I. Byrnes, eds, LNCIS vol. 286,
 Springer-Verlag, Berlin, 2002, pp.\,261--278.


\bibitem{StafMTNS02} O.J. Staffans, \m  Passive and conservative
 infinite-dimensional impedance and scattering systems (from a
 personal point of view), In {\em Mathematical Systems Theory in
 Biology, Communication, Computation, and Finance}, volume 134 of
 {\em IMA Volumes in Mathematics and its Applications}, pages
 375--414. Springer-Verlag, New York, 2002.

\bibitem{stafbook} O.J. Staffans, \m {\em Well-Posed Linear Systems},
 Cambridge University Press, Cambridge, UK, 2004. 

\bibitem{StWe} O.J. Staffans and G. Weiss, \m Transfer functions of
 regular linear systems, Part II: The system operator and the
 Lax-Phillips semigroup, \m {\it Trans. Amer. Math. Society}
 {\bf 354} (2002), pp.\,3329--3262.

\bibitem{SADG} V.L. Syrmos, C.T. Abdallah, P. Dorato and K.
 Grigoriadis, \m Static output feedback - A survey, \m {\it
 Automatica} {\bf 33} (1997), pp.\,125--137.

\bibitem{trigg2} R. Triggiani, \m Lack of uniform stabilization for
 noncontractive semigroups under compact perturbations, \m {\it
 Proc. Amer. Math. Soc.} {\bf 105} (1989) pp.\,375--383.

\bibitem{trigg1} R. Triggiani, \m Wave equation on a bounded domain
 with boundary dissipation: an operator approach, \m {\it J. Math.
 Anal. Appl.} {\bf 137} (1989), pp.\,438--461.

\bibitem{TuWe03} M. Tucsnak and G. Weiss, \m How to get a
 conservative well-posed linear system out of thin air. {Part II}:
 controllability and stability, \m {\em SIAM J. Control and Optim.}
 {\bf 42} (2003), pp.\,907--935. 
 
\bibitem{gang} J.A. Villegas, H. Zwart, Y. Le Gorrec, B. Maschke and
 A.J. van der Schaft, \m Stability and stabilization of a class of 
 boundary control systems, Proc. 44th IEEE Conference on Decision 
 and Control and the European Control Conference, Seville, Spain
 (2005), pp.\,3850--3855.

\bibitem{weiss1} G. Weiss, \m Admissibility of unbounded control
 operators, \m {\it SIAM J. Control and Optim.} {\bf 27} (1989),
 pp.\,527--545.


\bibitem{weiss10} G. Weiss, \m Transfer functions of regular linear
 systems, Part I: Characterizations of regularity, \m {\it Trans.
 Amer. Math. Society} {\bf 342} (1994), pp.\,827--854.

\bibitem{weiss12} G. Weiss, \m Regular linear systems with feedback, 
 \m {\it Mathematics of Control, Signals and Systems} \mm {\bf 7} 
 (1994), pp.\,23--57.

\bibitem{weiss63} G. Weiss, \m Optimal control of systems with a 
 unitary semigroup and with colocated control and observation, 
 \m {\it Systems \& Control Letters} {\bf 48} (2003), pp.\,329--340.

\bibitem{WCbeam} G. Weiss and R.F. Curtain, \m Exponential
 stabilization of a Rayleigh beam using colocated control, \m {\it
 IEEE Trans. on Automatic Control} {\bf 53} (2008), pp.\,643--654.

\bibitem{WeRe} G. Weiss and R. Rebarber, \m Optimizability and 
 estimatability for infinite-dimensional linear systems, \m {\it
 SIAM J. Control and Optim.} {\bf 39} (2001), pp.\,1204--1232.

\bibitem{WST} G. Weiss, O.J. Staffans and M. Tucsnak, \m Well-posed
 linear systems -a survey with emphasis on conservative systems, \m
 {\em Applied Mathematics and Computer Science} {\bf 11} (2001),
 pp.\,101--127.

\bibitem{WeTu} G. Weiss and M. Tucsnak, \m How to get a conservative
 well-posed linear system out of thin air. Part I: Well-posedness and
 energy balance, \m {\it ESAIM-COCV} {\bf 9} (2003), pp.\,247--274.

\bibitem{WeXu} G. Weiss and C-Z. Xu, \m Spectral properties of
 infinite-dimensional closed-loop systems, \m {\it Mathematics of 
 Control, Signals and Systems} {\bf 17} (2005), pp.\,153--172.

\bibitem{Will} J.C. Willems, \m Dissipative dynamical systems. Part I:
 General theory. Part II: Linear systems with quadratic supply rates,
 \m {\it Arch. Ration. Mech. Anal.} {\bf 45} (1972), pp.\,321--392.

\bibitem{You} Y. You, \m Dynamical boundary control of two-dimensional
 Petrovsky system: Vibrating rectangular plate. \m In A. Bensoussan
 and J.-L. Lions, editors, {\it Analysis and Optimization of Systems},
 volume 111 of LNCIS, pp.\,519--530. Springer-Verlag, Heidelberg, 
 1988.
 
\end{thebibliography}
\end{document}